# Two-Stage Data-Driven Contextual Robust Optimization: An End-to-End Learning Approach for Online Energy Applications

Carlos Gamboa[*,1,4], Alexandre Street[†,2,4], Davi Valladão[‡,2,4], and Bernardo Pagnocelli[§,3]

[1]Department of Electrical Engineering, PUC-Rio
[2]Department of Industrial Engineering, PUC-Rio
[3]SKEMA Business School
[4]Laboratory of Applied Mathematical Programming and Statistics (LAMPS), PUC-Rio

**Abstract**

Traditional end-to-end contextual robust optimization models are trained for specific contextual data, requiring complete retraining whenever new contextual information arrives. This limitation hampers their use in online decision-making problems such as energy scheduling, where multiperiod optimization must be solved every few minutes. In this paper, we propose a novel Data-Driven Contextual Uncertainty Set, which gives rise to a new End-to-End Data-Driven Contextual Robust Optimization model. For right-hand-side uncertainty, we introduce a reformulation scheme that enables the development of a variant of the Column-and-Constraint Generation (CCG) algorithm. This new CCG method explicitly considers the contextual vector within the cut expressions, allowing previously generated cuts to remain valid for new contexts, thereby significantly accelerating convergence in online applications. Numerical experiments on energy and reserve scheduling problems, based on classical test cases and large-scale networks (with more than 10,000 nodes), demonstrate that the proposed method reduces computation times and operational costs while capturing context-dependent risk structures. The proposed framework (model and method), therefore, offers a unified, scalable, and prescriptive approach to robust decision-making under uncertainty, effectively bridging data-driven learning and optimization.

## 1 Introduction

Large-scale infrastructure and network-constrained applications, such as energy and reserve scheduling in power systems, are increasingly impacted by geopolitical, environmental, and sustainability challenges. As nations seek energy independence by reducing reliance on external fuel sources, they drive the expansion of renewables, which inherently introduces uncertainty. Meanwhile, climate change exacerbates operational risks through more frequent and extreme weather events, which alter renewable generation patterns and strain existing infrastructure. A just and sustainable energy transition, therefore, requires systemic infrastructure and operational changes that balance climate mitigation efforts with geopolitical realities. In this evolving landscape, where long-term structural shifts combine with short-term operational variability, power system operators face unprecedented uncertainty, requiring robust and computationally efficient solutions for real-time decision-making.

[*]caagambaro@gmail.com
[†]street@ele.puc-rio.br
[‡]davimv@puc-rio.br
[§]bernardo.pagnoncelli@skema.edu

In an era of abundant data, it is natural for decision-makers to take into account *contextual information* before making decisions. By contextual information, we refer to the features observed together with a model's uncertain elements. For instance, in inventory problems, as demand is recorded, it is also useful to know the day of the week, the demand in the previous week, weather conditions, and macroeconomic variables at that point in time. While classical decision-making under uncertainty tools such as stochastic programming ([Shapiro et al., 2021]) and robust optimization ([Ben-Tal et al., 2009]) mostly ignored contextual information, we have been witnessing a significant growth in the number of publications that take them into account (see [Ban and Rudin, 2019, Bertsimas and Kallus, 2020, Elmachtoub and Grigas, 2022, Rahimian and Pagnoncelli, 2023, Estes and Richard, 2023]). In energy applications, some recent works have started to integrate contextual information within the decision-making model. In [Wahdany et al., 2023], the authors train neural networks, taking into account the system dispatch and their costs. [Casagrande et al., 2024] adopts a similar approach, integrating the MPC optimization as a differentiable optimization layer within the neural network. Contextual models gave rise to *end-to-end learning*, in which the models include the complete path from data to decisions. We refer the reader to the recent survey [Sadana et al., 2025] that describes the many different approaches in this emerging field.

Integrating state-of-the-art contextual models into operational decision-making problems in energy remains a significant challenge. These applications involve solving sequential network-constrained schedules that consider thousands of nodes and edges representing the transmission system limits as new information becomes available. Traditional approaches often require rerunning models entirely when uncertainty changes, leading to inefficiencies in online scheduling. Additionally, incorporating uncertainty into official models poses challenges due to the computational burden of complex stochastic and robust formulations. In this paper, we will focus on *data-driven robust optimization models*, which offer an attractive trade-off between uncertainty modeling and computational tractability.

The work of [Bertsimas and Brown, 2009] pioneered the construction of data-driven polyhedral uncertainty sets for uncertainty-constrained linear robust optimization problems. This approach leverages convex combinations of historical observations to define a family of generating probability measures for coherent risk measures, such as the Conditional Value at Risk (CVaR). Data-driven models, in general, have become central to quantifying uncertainty in both robust and distributionally robust optimization. By using finite samples from historical observations, these models encapsulate uncertainty either within specialized uncertainty sets or through discrete empirical probability distributions. This approach allows decision-makers to directly incorporate empirical data into the optimization process, thereby reducing reliance on stringent parametric assumptions and improving both the tractability and cost-effectiveness of the solutions. As a result, *data-driven robust optimization* has attracted significant attention and development in academic research, particularly in applications where uncertainty plays a pivotal role in operational decision-making.

The use of contextual information within robust optimization has the potential to reduce the conservatism of solutions by shrinking the uncertainty set given the available information. Contextual data-driven robust optimization models represent the state-of-the-art in robust optimization. These models can be classified into two categories: "predictive-then-optimize" and "end-to-end" models (see [Vanderschueren et al., 2022, Kotary et al., 2021] for discussions and comparisons).

The first category refers to models that use predictive procedures to construct a specialized uncertainty set that is conditioned to a given contextual information, which is then used to prescribe an optimal decision by solving the resulting robust optimization problem. For instance, [Sun et al., 2023] develops predictive-then-optimize algorithms for contextual linear programming optimization problems with objective uncertainty. [Patel et al., 2024] presents a model that delivers non-convex conditioned uncertainty sets with finite-sample guarantees, mitigating the conservatism of prescriptive decisions. [Peršak and Anjos, 2023] introduces a contextual data-driven distributionally robust optimization problem that utilizes supervised learning methods to construct the ambiguity set. [Chenreddy et al., 2022] develops a model that constructs the conditioned uncertainty set by using deep clustering methods to enhance the predictive estimation of the uncertainty and establishes a connection with Value-at-Risk optimization. Finally, in [Bertsimas et al., 2025], the authors presented a new method to construct uncertainty sets based on machine learning models and mixed-integer optimization. This approach leverages historical data on features and realizations of the uncertain parameters to deliver conditioned uncertainty set bounds based on the error minimization of parametric regression models. Notably, however, predictive-then-optimize methods deliver decisions that may disregard the existence of a better predictive model for that application (see [Dias Garcia et al., 2025] for further discussion in this theme).

As an alternative approach, researchers also explore the end-to-end category, which integrates the predictive estimation directly into the optimization procedure, enhancing the solution performance by accounting for the downstream application loss or revenue. This work focuses on this category of models. In [Chenreddy and Delage, 2024], an end-to-end algorithm for delivering contextual robust optimization solutions is presented; however, it considers only cost functions that are concave in the uncertainty. In [Yeh et al., 2024], an end-to-end model is introduced, assuming objective uncertainty to calibrate the uncertainty set using conformal prediction methods. This approach ensures a probabilistic guarantee—that is, the uncertain parameter belongs to the uncertainty set with a specified confidence level. The authors propose an exact differentiation method for the training procedure and explore different representations of the calibrated uncertainty set. To transform a prescriptive robust optimization model into an automated system that delivers decisions online based on contextual information, its adaptability to new observations is of paramount importance. Examples of such adaptability include [Fernandes et al., 2016], whose framework captures price dynamics and updates market patterns, and [Velloso et al., 2020], who addresses the complex (high-dimensional and non-Gaussian) spatial and temporal dependencies in the day-ahead unit commitment for power systems operation.

Interestingly, all previously reported end-to-end contextual robust models are contextually trained, i.e., are specifically trained for a given contextual information, thereby requiring complete retraining from scratch for each new incoming data. This imposes a significant computational burden and becomes a barrier for online applications such as energy scheduling. In such a problem, the sequential decision-making process relies on the optimization of a multiperiod security-constrained economic dispatch problem every 15-60 minutes (see [Ribeiro et al., 2023] for the granularities of each country). This becomes especially critical for the right-hand-side (RHS) uncertainty—generally associated with resource or demand uncertainty, such as renewable generation in energy scheduling problems. To the best of our knowledge, there is no previous work in the literature that develops an end-to-end data-driven robust optimization model that: 1) explicitly accounts for the contextual information vector and 2) allows previously trained models (i.e., previously obtained cuts) to be reused.

Our main objective is to propose a novel end-to-end framework for obtaining contextual robust data-driven decisions with recourse that can be applied to online decision-making. The approach can be viewed as an optimal predictive estimator that minimizes the final cost of context-dependent decisions, delivering not only contextually informed first-stage robust decisions but also cuts that can be reused with contextual information in subsequent periods. Covariate information is directly utilized to construct a contextually-adjusted data-driven polyhedral uncertainty set. This structure ensures that the resulting optimization problem remains computationally tractable while keeping worst-case scenarios grounded in data, thereby making them realistic [Velloso et al., 2020]. This problem can be efficiently solved using decomposition techniques. Based on a new reformulation of the second and third levels of the two-stage robust model, we develop a Column-and-Constraint Generation (CCG) method for problems with right-hand-side uncertainty that explicitly incorporates the contextual infor-mation vector. This salient feature of our methodology makes previously generated cuts also valid for future contextual information. Thus, the proposed method can start with previously generated valid cuts, quickly converging to the optimal solution for the updated contextual information, an essential feature for online optimization applications.

We test the performance of the proposed framework in an energy and reserve scheduling problem with uncertain renewable generation. To assess its performance on larger systems, we use different instances based on the well-known IEEE 118-bus test system, a 600-bus system created by combining two copies of the IEEE 300-bus system, and the 10,000-bus Grid Optimization Competition (GOC) test case described in [Babaeinejadsarookolaee et al., 2021]. The results demonstrate that the proposed reformulation and cutting sharing method significantly accelerates the CCG convergence. Additionally, our studies corroborate that incorporating contextual information reduces uncertainty, resulting in significant economic benefits. Specifically, the proposed framework achieves lower operating costs compared to those estimated by non-contextual robust optimization models. Our approach, therefore, stands out in its ability to offer contextually-informed robust solutions that are both practical and adaptable to changing conditions and uncertainties. This prescriptive analytics framework, grounded in the latest optimization techniques and data-driven insights, is poised to offer significant advantages in various fields where decision-making under uncertainty is critical.

Our key contributions are as follows:

1. We propose an innovative data-driven polyhedral uncertainty set, conditioned by contextual information, that enhances the adaptiveness of robust optimization models. This data-driven conditional uncertainty (DDCU) set can be interpreted as a data-driven nonparametric regression model embedded within the robust optimization problem. In this modeling approach, contextual information serves as regressors that locally adapt the geometry of an uncertainty set whose vertices are based on real data, extending the applicability of previous works [Velloso et al., 2020].

2. We present a novel end-to-end data-driven framework for Contextual Robust Optimization (CRO), giving rise to the Data-Driven Contextual Robust Optimization (DDCRO). This framework focuses on data-driven, two-stage robust optimization models with contextual information. In this framework, the need for a separate predictive step is eliminated through a unique modeling approach that leverages data-driven polyhedral uncertainty sets conditioned on contextual information.

3. For objective function uncertainty, we provide a convex reformulation of the DDCRO problem that efficiently scales with the number of historical observations used to construct the data-driven polyhedral uncertainty set, conditioned by contextual information.

4. For the right-hand-side uncertainty case, we propose a novel reformulation for the contextually adapted data-driven robust problem. This reformulation offers two key advantages for decomposition algorithms. Firstly, it results in lighter master problems and a tractable oracle for large-scale recourse problems. Secondly, our reformulation explicitly incorporates contextual information as variables within the cut expressions. This reformulation enhances efficiency in online applications, as cuts generated in previous runs remain valid and can be dynamically updated with new contextual information.

The structure of this paper is as follows. In Section 2, we describe the construction of the data-driven polyhedral uncertainty set, conditioned on contextual information. We also explore the concept of robustness of contextual information budget and provide guarantees for obtaining uncertainty-aware solutions. Section 3 introduces the DDCRO problem and presents reformulations for both objective and right-hand side uncertainty. The solution methodology for right-hand side uncertainty is detailed in Section 4, and extensions are presented in Section 5. In Section 6, we present an empirical study of the online energy application problem, along with the corresponding numerical results. Finally, conclusions are drawn in Section 7.

## 2 The proposed conditional uncertainty set

In robust optimization, uncertainty is modeled using a specially defined set, known as the uncertainty set, to find optimal solutions that are protected against worst-case scenarios for all possible realizations of the uncertainty within this set. In this work, we extend this concept to incorporate contextual information to forge the uncertainty set by exploring a novel technique for data-driven decision-making under uncertainty and contextual information known as Contextual Robust Optimization (CRO). This framework is designed to address robust optimization problems where the uncertainty set incorporates the most up-to-date side information derived from a set of covariates.

### 2.1 A new data-driven conditioned uncertainty (DDCU) set

We propose a polyhedral uncertainty set based on scenarios, following [Bertsimas and Brown, 2009, Fernandes et al., 2016, Bertsimas et al., 2018]. This uncertainty set is defined as the convex hull of a set of multivariate points that jointly represent scenarios of the covariates and the uncertainty. Among the various exogenous approaches that can be utilized for generating these scenarios, we suggest the implementation of a data-driven scheme. In this scheme, historical data of the covariates $\{\mathbf{x}_s\}_{s \in [S]}$ and uncertainty $\{\mathbf{y}_s\}_{s \in [S]}$ are used directly as scenarios This approach embeds pertinent information regarding the true underlying uncertainty process within each vertex of the uncertainty set, conditioned on the convex combination of the covariate scenarios being close enough to the most recent observation of contextual information $\mathbf{x}$. Mathematically, the proposed Data-Driven Conditioned Uncertainty (DDCU) set, denoted as $Y(\mathbf{x})$, is formulated as follows:

$$Y(\mathbf{x}) = \left\{ \tilde{\mathbf{y}} \in \mathbb{R}^{d_y} \;\middle|\; \begin{array}{l} \exists \boldsymbol{\theta} \in \mathbb{R}^S_+ : \quad \sum_{s \in [S]} \theta_s = 1 \\[4pt] \phantom{\exists \boldsymbol{\theta} \in \mathbb{R}^S_+ :} \quad \sum_{s \in [S]} \theta_s \mathbf{y}_s = \tilde{\mathbf{y}} \\[4pt] \phantom{\exists \boldsymbol{\theta} \in \mathbb{R}^S_+ :} \quad \sum_{s \in [S]} \theta_s \mathbf{x}_s = \mathbf{x} \end{array} \right\}. \qquad (1)$$

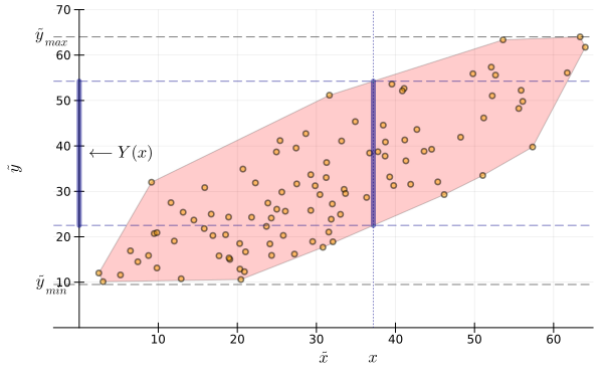

Figure 1: Illustrative Example: One-Dimensional Case

As illustrated in Figure 1, conditioning on contextual information allows us to reduce uncertainty more effectively compared to using an unconditioned data-driven polyhedral uncertainty set, as used in [Fernandes et al., 2016]. Without conditioning, the resulting robust optimization problem perceives the full span of uncertainty, which often leads to an overly conservative solution. As illustrated in Figure 1, once incorporated into the robust model, the proposed contextually adjusted uncertainty set (1) effectively embeds a nonparametric regression that links the uncertainties to the contextual information.

Moreover, the proposed DDCU set is adaptive to new observations of both covariates and uncertainty, similarly to its unconditioned counterpart. That is, historical realizations of covariates and un-certainty within a time span can be used as scenarios to construct the DDCU set. Then, in a rolling-window scheme, new observations of covariates and uncertainty can be incorporated into the dataset, while older observations are discarded. We highlight that the DDCU set defined in (1) represents an innovative approach to construct a data-driven conditioned uncertainty set within the robust optimization literature.

From definition (1), if $\mathbf{x}$ is not selected as a point that belongs to the convex hull of the covariate scenarios, then the set $Y(\mathbf{x})$ is empty. The emptiness of the uncertainty set is problematic, as it implies the absence of uncertainty quantification and renders the associated robust optimization problem infeasible. Moreover, definition (1) would be of limited practical use, particularly for high-dimensional problems. It is well known that, in high-dimensional spaces, a randomly selected point is, with high probability, located outside the convex hull of a fixed finite set of points [Vershynin, 2018][1].

To relax this somewhat stringent assumption, we extend our definition (1) to encompass a more general case in which the new contextual observation $\mathbf{x}$ does not necessarily belong to the convex hull of historical observations of the covariates. We do so by allowing points $\tilde{\mathbf{x}}$ within the convex hull of $\{\mathbf{x}_s\}_{s \in [S]}$ that are sufficiently close to the new contextual observation $\mathbf{x}$. We proposed the following relaxed definition, which supersedes definition (1) as the DDCU set:

$$Y_\Gamma(\mathbf{x}) = \left\{ \tilde{\mathbf{y}} \in \mathbb{R}^{d_y} \left| \begin{array}{l} \exists \boldsymbol{\theta} \in \mathbb{R}_+^S, \tilde{\mathbf{x}} \in \mathbb{R}^{d_x} : \quad \sum_{s \in [S]} \theta_s = 1 \\ \\ \sum_{s \in [S]} \theta_s \mathbf{y}_s = \tilde{\mathbf{y}} \\ \\ \sum_{s \in [S]} \theta_s \mathbf{x}_s = \tilde{\mathbf{x}} \\ \\ \|\tilde{\mathbf{x}} - \mathbf{x}\| \leq \Gamma \end{array} \right. \right\}. \tag{2}$$

In definition (2), a point $\mathbf{x}$ outside the convex hull of $\{\mathbf{x}_s\}_{s \in [S]}$ may generate a non-empty set $Y_\Gamma(\mathbf{x})$, as long as it is within $\Gamma$ of a point in the convex hull. Additionally, this formulation accounts for cases where there may be imprecision in the measurement of the observation $\mathbf{x}$.

The parameter $\Gamma$ can be interpreted as a contextual information robustness budget, which controls the level of robustness in the optimal decision, ranging from fully conditioned contextual information to fully unconditioned uncertainty robustness. Therefore, $\Gamma$ should be chosen small enough to ensure both the precision of the new observation of the contextual information $\mathbf{x}$ and the reduction of uncertainty that can be achieved by conditioning the set on $\mathbf{x}$.

Let us now discuss how different choices of $\Gamma$ affect the DDCU set. We start by defining $\Gamma_0$ as the distance between the new observation of contextual information $\mathbf{x}$ and the convex hull of the covariate data, represented as $\text{Conv}(\{\mathbf{x}_s\}_{s \in [S]})$. It is important to highlight that $\Gamma_0$ is the minimum value of the parameter $\Gamma$ required to ensure that $Y_\Gamma(\mathbf{x})$ is nonempty. Therefore, $\Gamma$ should be greater than or equal to $\Gamma_0$, that is, $\Gamma \geq \Gamma_0$. Obviously, if $\Gamma_0 = 0$, which implies that $\mathbf{x} \in \text{Conv}(\{\mathbf{x}_s\}_{s \in [S]})$, then $\Gamma$ should be set equal to $\Gamma_0$, as this ensures the maximum possible reduction of uncertainty for the contextual observation $\mathbf{x}$. When $\Gamma_0 = 0$, the only reason to choose $\Gamma > 0$ would be due to measurement imprecision of the contextual information $\mathbf{x}$.

In Figure 2a, we have an example of the case $\Gamma = \Gamma_0 = 0$. Geometrically, we have the intersection between the convex hull $\text{Conv}(\{(\mathbf{x}_s, \mathbf{y}_s)\}_{s \in [S]})$ of the data $\{(\mathbf{x}_s, \mathbf{y}_s)\}_{s \in [S]}$—illustrated by the blue polyhedron in the figure—and the line $\{(\tilde{\mathbf{x}}, \mathbf{y}) | \mathbf{y} \in \mathbb{R}^{d_y}\}$. The worst-case scenario $\tilde{\mathbf{y}}^*$ over $Y_\Gamma(\mathbf{x})$ is given by one of the endpoints of the resulting line segment from this intersection.

When $\mathbf{x}$ does not lie within $\text{Conv}(\{\mathbf{x}_s\}_{s \in [S]})$, we start by considering $\Gamma = \Gamma_0 > 0$. Following (2), we have that $\mathbf{x}$ is the center of the ball $B_1(\mathbf{x}, \Gamma)$ using the $\infty$-norm, illustrated by the green polyhedron in

---
[1] This is a well-known phenomenon referred to as the curse of dimensionality and the concentration of measure, which arises as a consequence of Milman's theorem. In high-dimensional spaces, the convex hull of a finite number of points occupies an exponentially small fraction of the ambient space's volume. Therefore, a new random point almost surely lies outside the convex hull.

Figure 2b. In this case, $\tilde{\mathbf{x}}$ may either coincide with the first coordinate of a vertex $(\tilde{\mathbf{x}}^*, \tilde{\mathbf{y}}^*)$ of the convex hull formed by the data $\{(\mathbf{x}_s, \mathbf{y}_s)\}_{s \in [S]}$ or lie on a facet of this polyhedron. Since it is highly probable that $(\tilde{\mathbf{x}}^*, \tilde{\mathbf{y}}^*)$ corresponds to a vertex of $\text{Conv}(\{(\mathbf{x}_s, \mathbf{y}_s)\}_{s \in [S]})$, the DDCU set is a singleton, denoted by $\{\tilde{\mathbf{y}}^*\}$.

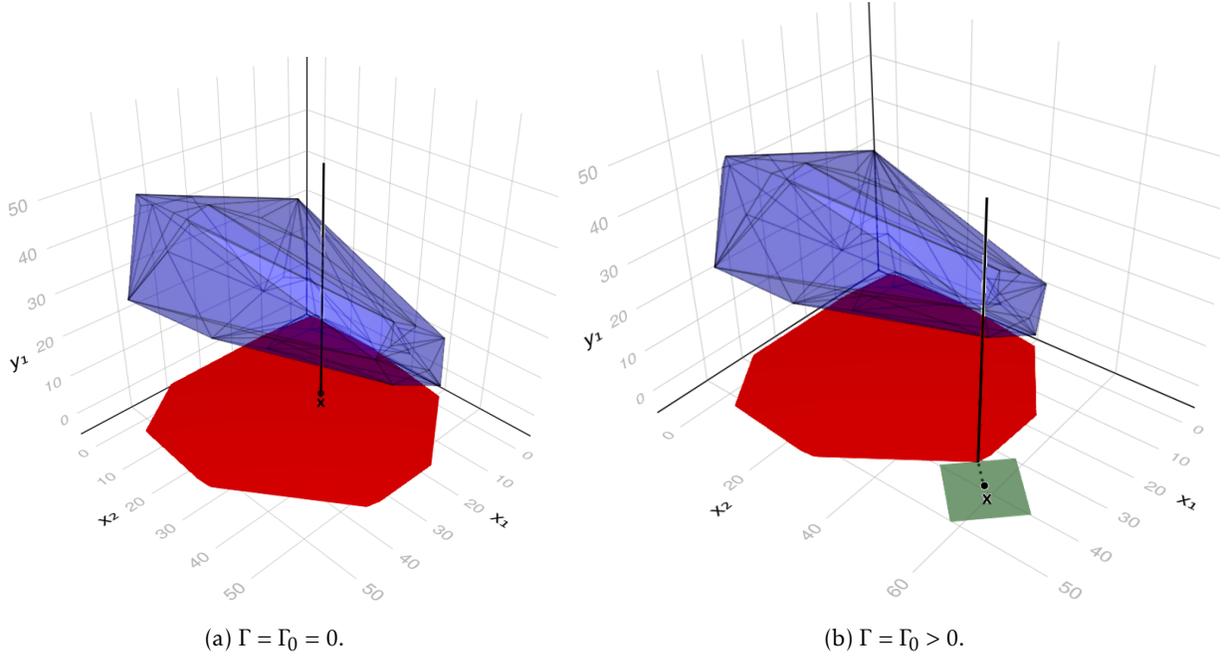

(a) $\Gamma = \Gamma_0 = 0$.

(b) $\Gamma = \Gamma_0 > 0$.

Figure 2: Case $\Gamma = \Gamma_0 \geq 0$.

For the case $\Gamma > \Gamma_0 > 0$, depicted in Figure 3, the DDCU set is defined as the projection onto the uncertainty dimension of the intersection between the polyhedron $\text{Conv}(\{(\mathbf{x}_s, \mathbf{y}_s)\}_{s \in [S]})$ and the cylinder $C = D \times \mathbb{R}^{d_y}$, whose base $D$ is given by $\text{Conv}(\{\mathbf{x}_s\}_{s \in [S]}) \cap B_1(\mathbf{x}, \Gamma)$—illustrated by the purple polyhedron in the figure. The worst-case scenario $\tilde{\mathbf{y}}^*$ corresponds to the second coordinate of a vertex $(\tilde{\mathbf{x}}^*, \tilde{\mathbf{y}}^*)$ of the restricted polyhedron $C \cap \text{Conv}(\{(\mathbf{x}_s, \mathbf{y}_s)\}_{s \in [S]})$ formed by this intersection.

It is worth mentioning that when the DDCU set $Y_\Gamma(\mathbf{x})$ results in a singleton $\mathbf{y}^*$—where uncertainty is quantified by a single scenario—the associated robust optimization problem becomes deterministic. In this situation, there is inherent risk because the optimal solution is not truly uncertainty-aware. In the following section, we will show conditions under which the set $Y_\Gamma(\mathbf{x})$ is not a singleton.

## 2.2 Robustness budget of contextual information

The parameter $\Gamma$ can be interpreted as a robustness budget of contextual information, i.e., a tolerance around the revealed contextual information $\mathbf{x}$ representing a certain value that captures imprecision or measurement errors in the observation of $\mathbf{x}$. Accordingly, $\Gamma$ should be calibrated to reduce uncertainty with respect to the discrepancy between the revealed contextual information and the observed realizations of the covariates in the data set $\{\mathbf{x}_s\}_{s \in [S]}$.

If $\Gamma$ is not sufficiently large, there may only exist one convex combination $\boldsymbol{\theta} = (\theta_1, \ldots, \theta_S)^\top$ such that $\sum_{s \in [S]} \theta_s \mathbf{x}_s = \tilde{\mathbf{x}}$, $\sum_{s \in [S]} \theta_s \mathbf{y}_s = \tilde{\mathbf{y}}$ and $\|\tilde{\mathbf{x}} - \mathbf{x}\| \leq \Gamma$. In other words, the uncertainty set $Y_\Gamma(\mathbf{x})$ consists of a singleton $\{\tilde{\mathbf{y}}\}$ or empty. The following result shows that $\Gamma$ can be chosen such that this does not happen.

**Proposition 1** *For a sufficiently large $\Gamma$, there exists a $\tilde{\mathbf{x}}$ such that the uncertainty set $Y_\Gamma(\mathbf{x})$ is not empty nor a singleton, implying the existence of uncertainty-aware solutions robust to contextual information.*

**Proof 2** *Proof of Proposition 1 For $\Gamma > \Gamma_0$, there exists a $\tilde{\mathbf{x}}$ interior point of the convex hull $\text{Conv}(\{\mathbf{x}_s\}_{s \in [S]})$ such that $\|\tilde{\mathbf{x}} - \mathbf{x}\| \leq \Gamma$. The set $\text{Conv}(\{\mathbf{x}_s\}_{s \in [S]})$ is the same as the projection of the convex hull $C = \text{Conv}(\{(\mathbf{x}_s, \mathbf{y}_s)\}_{s \in [S]}) \subseteq$*

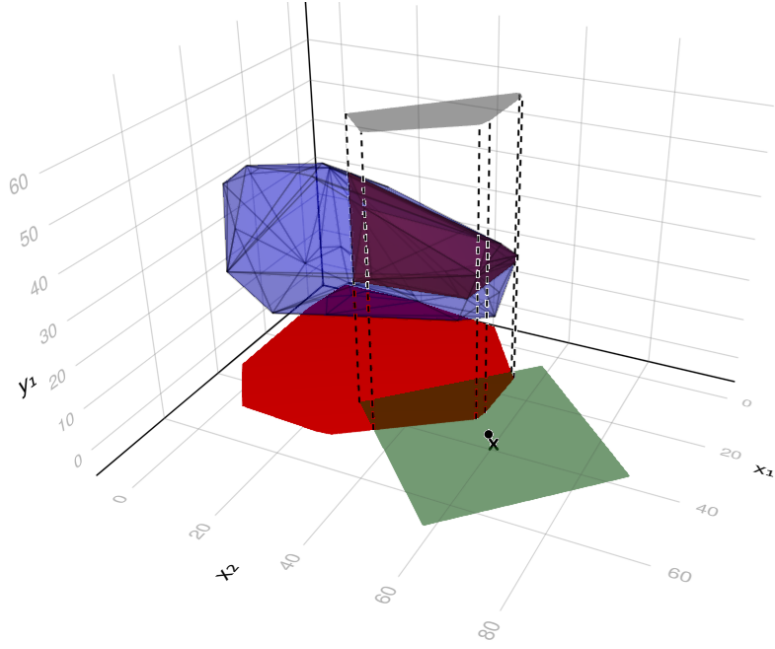

Figure 3: Case $\Gamma > \Gamma_0 > 0$.

$\mathbb{R}^{d_x+d_y}$ of the data $\{(\mathbf{x}_s, \mathbf{y}_s)\}_{s \in [S]}$ in the covariate space. By contradiction, suppose that $\tilde{\mathbf{y}}$ is the only $\mathbf{y}$-component that can pair with $\tilde{\mathbf{x}}$ to be in $C$. As $\tilde{\mathbf{x}}$ is an interior point of $Conv(\{\mathbf{x}_s\}_{s \in [S]})$, $(\tilde{\mathbf{x}}, \tilde{\mathbf{y}})$ is not a vertex of $C$ but it belongs to the boundary of $C$ since for any ball $B \subseteq \mathbb{R}^{d_x+d_y}$ centered at $(\tilde{\mathbf{x}}, \tilde{\mathbf{y}})$ we can find $\boldsymbol{\epsilon} \in \mathbb{R}_+^{d_y}$ such that the open interval $](\tilde{\mathbf{x}}, \tilde{\mathbf{y}} - \boldsymbol{\epsilon}), (\tilde{\mathbf{x}}, \tilde{\mathbf{y}} + \boldsymbol{\epsilon})[ \subseteq B$ and by hypothesis, $](\tilde{\mathbf{x}}, \tilde{\mathbf{y}} - \boldsymbol{\epsilon}), (\tilde{\mathbf{x}}, \tilde{\mathbf{y}} + \boldsymbol{\epsilon})[ \cap C = \{(\tilde{\mathbf{x}}, \tilde{\mathbf{y}})\}$. By the line segment principle, we can find $(\mathbf{x}', \mathbf{y}')$ and $(\mathbf{x}'', \mathbf{y}'')$ in the relative interior of $C$, $ri(C)$, with $\mathbf{x}' \neq \mathbf{x}''$ and $\mathbf{y}' \neq \mathbf{y}''$ such that the segment line $[(\mathbf{x}', \mathbf{y}'), (\mathbf{x}'', \mathbf{y}'')]$ lies on $ri(C)$. On the other hand, by hypothesis, the interior $int(C)$ of the set $C$ is not empty, which implies that $int(C) = ri(C)$. Therefore, $[(\mathbf{x}', \mathbf{y}'), (\mathbf{x}'', \mathbf{y}'')]$ lies in $int(C) = ri(C)$ and the point $(\tilde{\mathbf{x}}, \tilde{\mathbf{y}}') \in [(\mathbf{x}', \mathbf{y}'), (\mathbf{x}'', \mathbf{y}'')]$ for some $\tilde{\mathbf{y}}'$. So, if $\tilde{\mathbf{y}}' = \tilde{\mathbf{y}}$, the segment line $[(\mathbf{x}', \mathbf{y}'), (\mathbf{x}'', \mathbf{y}'')]$ does not lie on $int(C)$ because the point $(\tilde{\mathbf{x}}, \tilde{\mathbf{y}}') = (\tilde{\mathbf{x}}, \tilde{\mathbf{y}})$ is in the boundary of $C$, which is a contradiction. Otherwise, if $\tilde{\mathbf{y}}' \neq \tilde{\mathbf{y}}$ contradicts the fact that $\tilde{\mathbf{y}}$ is the only $\mathbf{y}$-component that can pair with $\tilde{\mathbf{x}}$ to be in $C$. All the possibilities of $\tilde{\mathbf{y}}'$ lead to a contradiction. Thus, the assumption that $\tilde{\mathbf{y}}$ is the only $\mathbf{y}$-component that can pair with $\tilde{\mathbf{x}}$ to be in $C$ where $\tilde{\mathbf{x}}$ is an interior point of the convex hull $Conv(\{\mathbf{x}_s\}_{s \in [S]})$ must be false. Therefore, there exists more than one $\tilde{\mathbf{y}} \in \mathbb{R}^{d_y}$ such that $(\tilde{\mathbf{x}}, \tilde{\mathbf{y}}) \in C = Conv(\{(\mathbf{x}_s, \mathbf{y}_s)\}_{s \in [S]})$ for $\tilde{\mathbf{x}}$ interior point of the convex hull $Conv(\{\mathbf{x}_s\}_{s \in [S]})$ such that $\|\tilde{\mathbf{x}} - \mathbf{x}\| \leq \Gamma$.

*Categorical features can be incorporated into the covariate vector, providing a mechanism to select scenarios for continuous features. Mathematically, the covariate vector $\tilde{\mathbf{x}}$, belonging to the convex hull $conv(\{\mathbf{x_s}\}_{s \in [S]})$ of the observed data (with $\mathbf{x}$ denoting a new covariate observation), can be written as*

$$\tilde{\mathbf{x}} := [\tilde{\mathbf{x}}^{cont}; \tilde{\mathbf{x}}^{cat}], \mathbf{x} := [\mathbf{x}^{cont}; \mathbf{x}^{cat}],$$

*where $\tilde{\mathbf{x}}^{cont}$ ($\mathbf{x}^{cont}$) represents the continuous feature components, and $\tilde{\mathbf{x}}^{cat}$ ($\mathbf{x}^{cat}$) represents the categorical feature components. For continuous features, the robustness budget associated with contextual information can vary between $\Gamma_0$ and an upper bound. In contrast, the robustness budget for categorical features must be zero, since categorical values in new observations cannot deviate from the observed ones. Thus, we require $\tilde{\mathbf{x}}^{cat} = \mathbf{x}^{cat}$.*

*It is possible, however, to assign different values of the parameter $\Gamma$ to the continuous and categorical components of $\tilde{\mathbf{x}}$ and $\mathbf{x}$. Moreover, if the norm $\|\cdot\|$ is the $\ell_1$-norm and a positive value of $\Gamma$ is (hypothetically) assigned to the categorical components, then the nearest integer approximation of $\Gamma$ specifies the maximum*

number of categorical features allowed to differ from the revealed contextual information in $\mathbf{x}^{cat}$. For simplicity of exposition, in what follows we will not distinguish between different values of $\Gamma$ for continuous and categorical components of the covariate vectors $\tilde{\mathbf{x}}$ and $\mathbf{x}$.

Since jointly considering categorical and continuous data induces a selection effect on the scenarios for continuous features, the data can be preprocessed to identify and select scenarios corresponding to each specific category. This approach leverages categorical contextual information in the covariate vector and may reduce the computational effort required to process the convex hull of the entire dataset $\{(\mathbf{x}_s, \mathbf{y}_s)\}_{s \in [S]}$.

# 3 Problem statement and proposed reformulations

We consider the following Data-Driven Contextual Robust Optimization (DDCRO) problem:

$$\min_{\mathbf{z} \in Z} \mathbf{c}^\top \mathbf{z} + \mathcal{Q}(\mathbf{z}, \mathbf{x}), \tag{3}$$

where $\mathbf{z}$ denotes the first-stage decision variable, $Z \subseteq \mathbb{R}^{d_z}$ denotes the feasible set of the first-stage problem, $\mathbf{c} \in \mathbb{R}^{d_z}$ represents the immediate cost, and the function $\mathcal{Q}(\mathbf{z}, \mathbf{x})$ provides the worst-case function value for any $\mathbf{z} \in Z$, conditioned on the new observation of covariates $\mathbf{x}$. Mathematically, the function $\mathcal{Q}(\mathbf{z}, \mathbf{x})$ is given by:

$$\mathcal{Q}(\mathbf{z}, \mathbf{x}) = \max_{\tilde{\mathbf{x}}, \tilde{\mathbf{y}} \in Y_\Gamma(\mathbf{x})} \min_{\mathbf{u} \geq 0} \mathbf{q}^\top \mathbf{u} \tag{4a}$$

$$s.t. \quad \mathbf{W}\mathbf{u} = \mathbf{h} - \mathbf{T}\mathbf{z}. \quad : \pi \tag{4b}$$

Problem (4) is a bilevel problem, in which the leader variables are $\tilde{\mathbf{x}}$ and $\tilde{\mathbf{y}}$, while the follower decides on $\mathbf{u}$. In the literature, the lower-level problem in (4) is known as the recourse or second-stage problem. Within this formulation, $\mathbf{u} \in \mathbb{R}^{d_u}$ is the second-stage decision variable, $\mathbf{q} \in \mathbb{R}^{d_u}$ represents the recourse cost, $\mathbf{T} \in \mathbb{R}^{d_h \times d_z}$ is the technology matrix, $\mathbf{W} \in \mathbb{R}^{d_h \times d_u}$ is the recourse matrix, $\mathbf{h} \in \mathbb{R}^{d_h}$ is the right-hand side vector, and $\pi$ is the dual variable of the constraint (4b). We use the notation $\tilde{\mathbf{y}} := (\mathbf{q}, \mathbf{h}, \mathbf{T}, \mathbf{W})$ to emphasize that the components of the random vector $\mathbf{y} \in \mathbb{R}^{d_y}$ can correspond to any component of the vectors $\mathbf{q} \in \mathbb{R}^{d_u}$, $\mathbf{h} \in \mathbb{R}^{d_h}$, or any entry of the matrices $\mathbf{T} \in \mathbb{R}^{d_h \times d_z}$ and $\mathbf{W} \in \mathbb{R}^{d_h \times d_u}$.

To find a solution robust to contextual information, we first focus on the problem (4). This problem is over the set $Y_\Gamma(\mathbf{x})$ for realizations $\tilde{\mathbf{y}}$ of the uncertainty, for which there exists a convex combination of data points $\{\mathbf{x}_s\}_{s \in [S]}$ that is equal to a vector $\tilde{\mathbf{x}} \in \mathbb{R}^{d_x}$ that is within a distance $\Gamma$ of $\mathbf{x}$. It is essential to emphasize that the point $\tilde{\mathbf{x}}$ is a decision variable. Consequently, if $\Gamma$ is too large, we would end up with the convex hull of the original data set, which would render the contextual information irrelevant. On the other hand, as we discussed in Section 2, if $\Gamma$ is not sufficiently large, problem (4) may be infeasible, as there might not be a point $\tilde{\mathbf{x}}$ in the convex hull of the data $\{\mathbf{x}_s\}_{s \in [S]}$ in the covariate space whose distance to $\mathbf{x}$ is smaller than or equal to $\Gamma$.

From this point forward, we use the contextual information robustness budget as a percentage $\delta$ of the distance from the new observation of covariates $\mathbf{x}$ to the convex hull defined by the data $\{\mathbf{x}_s\}_{s \in [S]}$, represented by $\Gamma_0$, i.e., $\Gamma = (1 + \delta)\Gamma_0$. Moreover, we will focus on the righ-hand-side uncertainty case, i.e., $\tilde{\mathbf{y}} := (\mathbf{h}, \mathbf{T})$, where $\mathbf{h} = \sum_{s \in [S]} \theta_s \mathbf{h}_s$ and $\mathbf{T} = \sum_{s \in [S]} \theta_s \mathbf{T}_s$ for some $\boldsymbol{\theta} = (\theta_1, \ldots, \theta_S)^\top \in \mathbb{R}_+^S$ such that $\sum_{s \in [S]} \theta_s = 1$ and $\sum_{s \in [S]} \theta_s \mathbf{x}_s = \mathbf{x}$. Let us rewrite problem (4) for the subsequent sections in the following equivalent form, hereinafter referred to as **P-Bilevel**[2]:

$$\mathcal{Q}(\mathbf{z}, \mathbf{x}) = \max_{\varepsilon, \boldsymbol{\theta} \geq 0} \min_{\mathbf{u} \geq 0} \left\{ \mathbf{q}^\top \mathbf{u} \;\middle|\; \mathbf{W}\mathbf{u} = \sum_{s \in [S]} \theta_s \big(\mathbf{h}_s - \mathbf{T}_s \mathbf{z}\big) \right\} \tag{5a}$$

$$s.t. \quad \sum_{s \in [S]} \theta_s = 1 \qquad : \rho, \tag{5b}$$

---

[2]We use the notation P-Bilevel to refer to the formulation (5) of the worst-case value function $\mathcal{Q}(\mathbf{z}, \mathbf{x})$, following the primal formulation of the lower-level problem in the bilevel optimization model that defines this function.

$$\sum_{s \in [S]} \theta_s \mathbf{x}_s - \varepsilon = \mathbf{x} \qquad\qquad : \gamma, \qquad (5c)$$

$$\|\varepsilon\| \le \Gamma, \qquad (5d)$$

where $\varepsilon = \tilde{\mathbf{x}} - \mathbf{x}$, and $\rho \in \mathbb{R}$ and $\gamma \in \mathbb{R}^{d_x}$ are the dual variables for constraint (5b) and (5c), respectively. We emphasize that the P-Bilevel formulation is a novel way to consider contextual information into two-stage data-driven robust models, giving rise to our second contribution, the DDCRO model.

## 3.1 Recourse problem with right-hand-side uncertainty

In the case of right-hand-side uncertainty, we derive an alternative and tractable bilevel reformulation of problem (5). A salient feature of this new formulation is that the contextual information vector $\mathbf{x}$ is explicitly included in the recourse objective function. This fact will unlock the possibility of sharing previously generated cuts for different contextual information, which will be critical for accelerating the CCG algorithm in online-optimization applications. Additionally, we propose an algorithmic scheme that enables the solution of the original problem (3) using a decomposition method. The bilevel problem (5) can be particularly challenging due to the complexity of the lower-level problem, especially in large-scale recourse settings. Drawing from [Saif and Delage, 2021], we present the following theorem:

**Theorem 3** *The P-Bilevel formulation (5) is equivalent to the following problem, hereinafter referred to as D-Bilevel[3]:*

$$\mathcal{Q}(\mathbf{z}, \mathbf{x}) = \max_{\pi} \quad \min_{\rho, \gamma} \left\{ \rho + \mathbf{x}^\top \gamma + \Gamma \|\gamma\|_* \,\middle|\, \rho + \mathbf{x}_s^\top \gamma \ge \pi^\top (\mathbf{h}_s - \mathbf{T}_s \mathbf{z}), \forall s \in [S] \right\}$$
$$\text{s.t.} \quad \mathbf{W}^\top \pi \le \mathbf{q}, \qquad (6)$$

*where $\|\cdot\|_*$ denotes the dual norm of the norm $\|\cdot\|$.*

**Proof 4** *Proof of Theorem (3)*

*By using the dual formulation of the lower (inner) problem in (5), we establish the following equivalence:*

$$\mathcal{Q}(\mathbf{z}, \mathbf{x}) = \max_{\varepsilon, \theta \ge 0} \left\{ \max_{\pi} \left\{ \sum_{s \in [S]} \theta_s \pi^\top (\mathbf{h}_s - \mathbf{T}_s \mathbf{z}) \,\middle|\, \mathbf{W}^\top \pi \le \mathbf{q} \right\} \,\middle|\, \begin{array}{l} \sum_{s \in [S]} \theta_s = 1 \\ \sum_{s \in [S]} \theta_s \mathbf{x}_s - \varepsilon = \mathbf{x} \\ \|\varepsilon\| \le \Gamma \end{array} \right\}. \qquad (7)$$

*Since the decision variables in the lower-level problem do not affect the constraints and the objective function of the upper-level problem (hierarchical independence), we can interchange the maximizations as follows:*

$$\max_{\pi} \left\{ \max_{\varepsilon, \theta \ge 0} \left\{ \sum_{s \in [S]} \theta_s \pi^\top (\mathbf{h}_s - \mathbf{T}_s \mathbf{z}) \,\middle|\, \begin{array}{ll} \sum_{s \in [S]} \theta_s = 1 & : \rho \\ \sum_{s \in [S]} \theta_s \mathbf{x}_s - \varepsilon = \mathbf{x} & : \gamma \\ \|\varepsilon\| \le \Gamma & \end{array} \right\} \,\middle|\, \mathbf{W}^\top \pi \le \mathbf{q} \right\}. \qquad (8)$$

*Note that in (8), we encounter nonlinearity due to the product of $\theta_s \pi$, $s \in [S]$. To mitigate this nonlinearity, we dualize the lower-level problem in (8), leading to the following bilevel formulation:*

$$\mathcal{Q}(\mathbf{z}, \mathbf{x}) = \max_{\pi} \quad \min_{\rho, \gamma} \left\{ \rho + \mathbf{x}^\top \gamma + \Gamma \|\gamma\|_* \,\middle|\, \rho + \mathbf{x}_s^\top \gamma \ge \pi^\top (\mathbf{h}_s - \mathbf{T}_s \mathbf{z}), \forall s \in [S] \right\}$$
$$\text{s.t.} \quad \mathbf{W}^\top \pi \le \mathbf{q}, \qquad (9)$$

*which is what we set out to prove.*

---
[3] We use the notation D-Bilevel to refer to the formulation (6) of the worst-case value function $\mathcal{Q}(\mathbf{z}, \mathbf{x})$, inspired by the incorporation of dual information—such as the dual polyhedron $\{\pi \in \mathbb{R}^{d_h} \mid \mathbf{W}^\top \pi \le \mathbf{q}\}$—into the bilevel optimization model that defines this function

*In the cases of using norm sup or norm 1 in (8), Proposition 3 allows us to work with a linear bilevel problem whose lower problem is characterized by simple constraints defined solely by the data $\{(\mathbf{x}_s, \mathbf{h}_s, \mathbf{T}_s)\}_{s \in [S]}$. Notably, the feasible set remains independent of contextual information, while the objective function explicitly incorporates contextual information $\mathbf{x}$.*

*By integrating the lower-level KKT conditions into the upper-level problem, as discussed in [Saif and Delage, 2021, Garcia et al., 2023], problem (9) can be restructured into a mixed integer linear program (MILP) as follows:*

$$\max_{\substack{\mathbf{B}, \rho, \gamma, \pi, \\ \theta \geq 0}} \quad \rho + \mathbf{x}^\top \gamma + \Gamma \|\gamma\|_* \tag{10a}$$

$$\text{s.t.} \quad \sum_{s \in [S]} \theta_s = 1, \tag{10b}$$

$$\sum_{s \in [S]} \theta_s \mathbf{x}_s = \mathbf{x}, \tag{10c}$$

$$\rho + \mathbf{x}_s^\top \gamma \geq \pi^\top \left( \mathbf{h}_s - \mathbf{T}_s \mathbf{z} \right), \quad \forall s \in [S], \tag{10d}$$

$$0 \leq \rho + \mathbf{x}_s^\top \gamma - \pi^\top \left( \mathbf{h}_s - \mathbf{T}_s \mathbf{z} \right) \leq M B_s, \quad \forall s \in [S], \tag{10e}$$

$$0 \leq \theta_s \leq M(1 - B_s), \quad \forall s \in [S], \tag{10f}$$

$$\mathbf{W}^\top \pi \leq \mathbf{q}, \tag{10g}$$

$$B_s \in \{0, 1\}, \quad \forall s \in [S], \tag{10h}$$

$$\rho \in \mathbb{R}, \gamma \in \mathbb{R}^{d_x}, \pi \in \mathbb{R}^{m_u}, \theta \in \mathbb{R}_+^S, \mathbf{B} \in \mathbb{R}^S, \tag{10i}$$

*where $\mathbf{B} = (B_1, \ldots, B_S)^\top$ are binary variables, $M$ is a constant, constraints (10b)-(10c) are dual feasibility, constraint (10d) is primal feasibility, and the disjunctive constraints (10e)-(10f) represents the complementary slackness condition of the lower-level problem in (9).*

*In bilevel optimization, the main challenge typically arises from managing the complementary slackness condition. However, given the interchange and duality relationship between (8) and (9), addressing the complementary slackness condition of the lower-level problem in (9) tends to be easier than handling it in the P-Bilevel formulation or for the dual polyhedron as in (7). This is due to the reduced number of constraints and lower dimensionality in the former cases. Consequently, our new reformulation (9) is particularly suited for complex recourse problems requiring optimization with the latest contextual information. However, it is important to note that the literature, as reviewed in [Garcia et al., 2023], offers a broad spectrum of solution techniques for the bilevel problem represented by (9). While valuable, an in-depth discussion of these methodologies is beyond the scope of this work.*

*In particular, when the uncertainty is only in the vector $\mathbf{h}$, since $\rho$ is unrestricted, we have that formulation (9) reduces to the following problem:*

$$\mathcal{Q}(\mathbf{z}, \mathbf{x}) = \max_{\pi} \quad -\pi^\top \mathbf{T} \mathbf{z} + \min_{\rho, \gamma} \left\{ \rho + \mathbf{x}^\top \gamma + \Gamma \|\gamma\|_* \,\middle|\, \rho + \mathbf{x}_s^\top \gamma \geq \pi^\top \mathbf{h}_s, \forall s \in [S] \right\}$$
$$\text{s.t.} \quad \mathbf{W}^\top \pi \leq \mathbf{q}. \tag{11}$$

*In this case, the decision variable $\mathbf{z} \in Z$ does not impact the feasible set of the lower-level problem as defined in (11). Its influence is limited solely to the objective function of the upper-level problem. From this point forward, the newly restructured formulation (11) of the bilevel problem (4) when uncertainty is confined to $\mathbf{h}$, will be denoted as $D_h$-Bilevel. This reformulation uniquely utilizes dual information in both the lower- and upper-level problems, offering a distinct approach to addressing the complexities inherent in bilevel optimization.*

# 4 Solution methodology for right-hand side uncertainty

*In this section, we present the proposed numerical schemes to solve the two-stage conditional robust optimization problem (3) for the right-hand side uncertainty case, i.e., $\tilde{\mathbf{y}} := (\mathbf{h}, \mathbf{T})$. We show that this problem can be addressed using the classical CCG method as well as a novel approach, the Column-and-Constraint Generation method with explicit contextual information (Contextual CCG). We first present the classical method and then describe our novel solution methodology. We first present the classical method and then describe our novel solution methodology. In the notation, symbols with subscripts "k" and superscripts "(k)" indicate new variables and variable values corresponding to the k-th scenario selected by the solution method.*

*Let us rewrite the problem (3) as follows:*

$$\min_{\mathbf{z},\alpha} \mathbf{c}^\top \mathbf{z} + \alpha \tag{12a}$$

$$\text{s.t.} \quad \alpha \geq \mathcal{Q}(\mathbf{z}, \mathbf{x}) \tag{12b}$$

$$\mathbf{z} \in Z. \tag{12c}$$

*With the sets, $Y_\Gamma(\mathbf{x})$ and $\Pi := \{\boldsymbol{\pi} \in \mathbb{R}^{d_h} \mid \mathbf{W}^\top \boldsymbol{\pi} \leq \mathbf{q}\}$ we can rewrite the worst-case function value $\mathcal{Q}(\mathbf{z}, \mathbf{x})$ based on the two possible formulations (P-Bilevel or D-Bilevel) for the second stage. For instance, according to the P-Bilevel formulation,*

$$\mathcal{Q}(\mathbf{z},\mathbf{x}) = \max_{\mathbf{y}} \left\{ Q^P(\mathbf{z},\mathbf{y}) \mid \mathbf{y} \in Y_\Gamma(\mathbf{x}) \right\}.$$

*Equivalently, according to the D-Bilevel formulation,*

$$\mathcal{Q}(\mathbf{z},\mathbf{x}) = \max_{\boldsymbol{\pi}} \left\{ Q^D(\mathbf{z},\boldsymbol{\pi},\mathbf{x}) \mid \boldsymbol{\pi} \in \Pi \right\},$$

*where*

$$Q^P(\mathbf{z},\mathbf{y}) := \min_{\mathbf{u} \geq 0} \left\{ \mathbf{q}^\top \mathbf{u} \,\middle|\, \mathbf{W}\mathbf{u} = \mathbf{h} - \mathbf{T}\mathbf{z} \right\},$$

*and*

$$Q^D(\mathbf{z},\boldsymbol{\pi},\mathbf{x}) := -\boldsymbol{\pi}^\top \mathbf{T}\mathbf{z} + \min_{\rho,\gamma} \left\{ \rho + \mathbf{x}^\top \boldsymbol{\gamma} + \Gamma \|\boldsymbol{\gamma}\|_* \,\middle|\, \rho + \mathbf{x}_s^\top \boldsymbol{\gamma} \geq \boldsymbol{\pi}^\top \mathbf{h}_s, \forall s \in [S] \right\}.$$

*As described in [Gamboa et al., 2021], a decomposition method can be implemented using a master-oracle scheme, as seen in [Velloso et al., 2020] for an application of this technique to a two-stage robust Unit Commitment problem. In our context, the master problem is a relaxation of the problem (12)*

*To derive a valid relaxation of the original problem (12), we can consider either a subset*

$$\{\mathbf{y}^{(k)} := (\mathbf{h}^{(k)}, \mathbf{T}^{(k)})\}_{k \in [K]} \subseteq Y_\Gamma(\mathbf{x})$$

*or*

$$\{\boldsymbol{\pi}^{(k)}\} \subseteq \Pi,$$

*defined at the iteration K of the iterative procedure. Once the solution to the master (relaxed) problem is obtained, the oracle identifies the worst infeasibility on the uncertainty set $Y_\Gamma(\mathbf{x})$ conditioned on the new observation $\mathbf{x}$ of the covariates or an extreme point of the polyhedron $\Pi$. It then adds the corresponding constraint (or block of constraints and variables) to the master problem. This iterative procedure stops when the oracle determines that the master solution is feasible for the original problem (3).*

*Note that, unlike $Q^P$, $Q^D(\mathbf{z}, \boldsymbol{\pi}, \mathbf{x})$ explicitly depends on $\mathbf{x}$ in the objective function, which is a salient feature of this formulation. This feature unlocks our second contribution: the possibility of using valid cuts obtained for different values of $\mathbf{x}$ within a decomposition scheme, thereby accelerating the convergence of the algorithm in online applications. This important aspect of our work is explored in the next section.*

## 4.1 Master-oracle method

We propose a solution methodology to address the problem (12) using a master-oracle numerical scheme. We develop an iterative procedure based on lower and upper bounding approximations of the problem (12), which converges to its optimal value and optimal solution.

Considering the subset $\{\mathbf{y}^{(k)} := (\mathbf{h}^{(k)}, \mathbf{T}^{(k)})\}_{k \in [K]} \subseteq Y_\Gamma(\mathbf{x})$ defined at the iteration $K$ of the iterative procedure and the associated subset $\{\boldsymbol{\pi}^{(k)}\}_{k \in [K]} \subseteq \Pi$, where $\boldsymbol{\pi}^{(k)} \in \arg\max_\pi \left\{ \boldsymbol{\pi}^\top (\mathbf{h} - \mathbf{T}\mathbf{z}) \middle| \boldsymbol{\pi} \in \Pi \right\}$ for all $k \in [K]$, we obtain the following problems:

$$\min_{\mathbf{z}, \alpha} \mathbf{c}^\top \mathbf{z} + \alpha \tag{13a}$$

$$\text{s.t.} \quad \alpha \geq Q^P(\mathbf{z}, \mathbf{y}^{(k)}), \quad \forall k \in [K] \tag{13b}$$

$$\mathbf{z} \in Z, \tag{13c}$$

and

$$\min_{\mathbf{z}, \alpha} \mathbf{c}^\top \mathbf{z} + \alpha \tag{14a}$$

$$\text{s.t.} \quad \alpha \geq Q^D(\mathbf{z}, \boldsymbol{\pi}^{(k)}, \mathbf{x}), \quad \forall k \in [K] \tag{14b}$$

$$\mathbf{z} \in Z, \tag{14c}$$

which represent relaxations of (12). Therefore, its optimal value is a valid lower bound, LB, for the optimal value of problem (12). Henceforth, we call the relaxed problem (13) or (14) as the master problem.

For the current optimal solution $(\mathbf{z}^{(K)}, \alpha^{(K)})$ of the master problem (13) or (14), we use the oracle problem $Q(\mathbf{z}^{(K)}, \mathbf{x})$—as defined by (4)—to find the worst-case uncertainty realization $\mathbf{y}^{(K)} \in Y_\Gamma(\mathbf{x})$ and its associated $\boldsymbol{\pi}^{(k)} \in \Pi$ among the restricted feasible set, which gets enlarged at each iteration. We summarize the master-oracle algorithm scheme in the following pseudo-code:

---

**Algorithm 1** Master-oracle algorithm

1: **Initialization**: Set $K = 0$, $UB^{(K)} \longleftarrow +\infty$ and $LB^{(K)} \longleftarrow -\infty$,
2: **while** $UB^{(K)} - LB^{(K)} > \text{tol}$ **do**
3:     **if** K = 0 **then**
4:         Solve the master problem (13) or (14) without $\alpha$ and the constraint (13b) or (14b);
5:         Store the master solution: $(LB^{(K)}, \mathbf{z}^{(K)})$;
6:         Set $\alpha^{(K)} = -\infty$
7:     **else**
8:         Solve the master problem (13) or (14);
9:         Store the master solution: $(LB^{(K)}, \mathbf{z}^{(K)}, \alpha^{(K)})$;
10:    **end if**
11:    Solve the oracle problem (4) for $\mathbf{z}^{(K)}$, i.e., find the optimal value $Q(\mathbf{z}^{(K)}, \mathbf{x})$;
12:    Store the oracle solution: $(\boldsymbol{\pi}^{(K)}, \mathbf{h}^{(K)}, \mathbf{T}^{(K)}) \in \{(\boldsymbol{\pi}^{(k)}, \mathbf{h}^{(k)}, \mathbf{T}^{(k)})\}_{k \in [K-1]}$;
13:    Set $UB \longleftarrow \max\{Q(\mathbf{z}^{(K)}, \mathbf{x}), \alpha^{(K)}\}$
14:    Update $K \longleftarrow K + 1$, add the constraint (13b)
15: **end while**
16: **return** $\mathbf{z}^{(K)}$, $UB^{(K)}$, $LB^{(K)}$

---

### 4.1.1 Classical Column-and-Constraint Generation method

Let $\boldsymbol{\theta}^{(k)} = (\theta_s^{(k)})_{s \in [S]}$ be the optimal solution corresponding to $\mathbf{z}^{(k)}$, obtained from the upper-level problem in the P-Bilevel formulation. Alternatively, this can be viewed as the dual variable for the constraint $\rho + \mathbf{x}_s^\top \boldsymbol{\gamma} \geq \boldsymbol{\pi}^\top (\mathbf{h}_s - \mathbf{T}_s \mathbf{z}), s \in [S]$, in the lower-level problem of the D-Bilevel formulation. Such a solution facilitates the *integration of primal information from the lower problem in (4) with endogenous contextual information into the master problem, ensuring feasibility. This integration is a fundamental aspect of the CCG algorithm.*

To define the master problem in the CCG algorithm, let $\mathbf{h}^{(k)} = \sum_{s \in [S]} \theta_s^{(k)} \mathbf{h}_s$ and $\mathbf{T}^{(k)} = \sum_{s \in [S]} \theta_s^{(k)} \mathbf{T}_s$. The master problem is then formulated as follows:

$$\min_{\mathbf{z}, \mathbf{u}_k, \alpha} \mathbf{c}^\top \mathbf{z} + \alpha \tag{15a}$$

$$\text{s.t.} \quad \alpha \geq \mathbf{q}^\top \mathbf{u}_k, \quad \forall k \in [K], \tag{15b}$$

$$\mathbf{W}\mathbf{u}_k = \mathbf{h}^{(k)} - \mathbf{T}^{(k)} \mathbf{z}, \quad \forall k \in [K], \tag{15c}$$

$$\mathbf{u}_k \geq 0, \quad \forall k \in [K], \tag{15d}$$

$$\mathbf{z} \in Z. \tag{15e}$$

*It is important to note that with each iteration of the CCG algorithm, the computational complexity to solve the master problem tends to increase. However, there is typically a balance to be struck between the computational load and the number of iterations required to ensure feasibility. Within this framework, Algorithm 1 must be executed by replacing the master problem (13) with its counterpart (15) at line 8.*

### 4.1.2 Contextual Column-and-Constraint Generation method

*In an alternative approach, the optimal solution $\boldsymbol{\pi}^{(k)}$, corresponding to the sub-optimal solution $\mathbf{z}^{(k)}$, enables us to approximate the worst-case value function $\mathcal{Q}(\mathbf{z}, \mathbf{x})$. This is achieved by utilizing the lower-level problem in the D-Bilevel formulation to ensure feasibility. The approximation is characterized by the use of cuts that explicitly exhibit contextual information. Equipped with the solution $\boldsymbol{\pi}^{(k)}$, we are able to accurately approximate the worst-case value function $\mathcal{Q}(\mathbf{z}, \mathbf{x})$. This involves generating new variables $\rho_k$ and $\boldsymbol{\gamma}_k$ and incorporating a set of linear constraints, $\rho_k + \mathbf{x}_s^\top \boldsymbol{\gamma}_k \geq (\boldsymbol{\pi}^{(K)})^\top (\mathbf{h}_s - \mathbf{T}_s \mathbf{z})$, for all $s \in [S]$ and for each iteration $k \in [K]$. These constraints are added to the master problem to modify the feasibility of the current sub-optimal solution $\mathbf{z}^{(K)}$ with each iteration k of the algorithm. Hereafter, this algorithm will be known as the Contextual CCG. The master problem for the Contextual CCG algorithm is defined as follows:*

$$\min_{\mathbf{z}, \alpha, \rho_k, \boldsymbol{\gamma}_k} \mathbf{c}^\top \mathbf{z} + \alpha \tag{16a}$$

$$\text{s.t.} \quad \alpha \geq \rho_k + \mathbf{x}^\top \boldsymbol{\gamma}_k + \Gamma \|\boldsymbol{\gamma}_k\|_*, \quad \forall k \in [K], \tag{16b}$$

$$\rho_k + \mathbf{x}_s^\top \boldsymbol{\gamma}_k \geq (\boldsymbol{\pi}^{(k)})^\top (\mathbf{h}_s - \mathbf{T}_s \mathbf{z}), \quad \forall s \in [S], \forall k \in [K], \tag{16c}$$

$$\mathbf{z} \in Z. \tag{16d}$$

*The ability to capture the worst-case value function $\mathcal{Q}(\mathbf{z}, \mathbf{x})$ using cuts infused with explicit contextual information from the lower-level problem in the D-Bilevel formulation renders the Contextual CCG algorithm highly adaptive. This adaptability implies that the algorithm can be initialized using the converged master problem, tailored to a previously observed contextual information, for a new observation of the covariates. This significant characteristic is encapsulated in the following proposition:*

**Proposition 5** *Fix an index set of columns $[K] = \{1, \ldots, K\}$ and the extreme points $\boldsymbol{\pi}^{(k)} \in \Pi$, $k \in [K]$, that have been obtained from previous oracle calls in the CCG algorithm. For any new observation of the covariates $\mathbf{x} \in \mathbb{R}^d$ and for every decision vector $\mathbf{z} \in Z$, the inequalities*

$$\alpha \geq \rho_k + \mathbf{x}^\top \boldsymbol{\gamma}_k + \Gamma \|\boldsymbol{\gamma}_k\|_* \quad \forall k \in [K], \tag{17a}$$

$$\rho_k + \mathbf{x}_s^\top \boldsymbol{\gamma}_k \geq (\boldsymbol{\pi}^{(k)})^\top (\mathbf{h}_s - \mathbf{T}_s \mathbf{z}) \quad \forall s \in [S], k \in [K] \tag{17b}$$

*are valid for the original upper-level problem*

$$\min_{\mathbf{z},\alpha} \quad \mathbf{c}^\top \mathbf{z} + \alpha \tag{18a}$$

$$\text{s.t.} \quad \alpha \geq \mathcal{Q}(\mathbf{z},\mathbf{x}), \tag{18b}$$

$$\mathbf{z} \in Z. \tag{18c}$$

**Proof 6** *Proof of Proposition 5 Fix an arbitrary $k \in [K]$ and $\mathbf{z} \in Z$. Define the inner optimisation problem of (17b) for the current covariate vector $\mathbf{x}$ and the stored column $\boldsymbol{\pi}^{(k)}$:*

$$(\rho_k^*, \boldsymbol{\gamma}_k^*) := \min_{\rho,\boldsymbol{\gamma}} \left\{ \rho + \mathbf{x}^\top \boldsymbol{\gamma} + \Gamma \|\boldsymbol{\gamma}\|_* \mid \rho + \mathbf{x}_s^\top \boldsymbol{\gamma} \geq (\boldsymbol{\pi}^{(k)})^\top (\mathbf{h}_s - \mathbf{T}_s \mathbf{z}), \ \forall s \in [S] \right\}.$$

*By construction the pair $(\rho_k^*, \boldsymbol{\gamma}_k^*)$ satisfies the constraints (17b). Set $(\rho_k, \boldsymbol{\gamma}_k) := (\rho_k^*, \boldsymbol{\gamma}_k^*)$ in (17) and denote the resulting objective value by*

$$v_k(\mathbf{z},\mathbf{x}) := \rho_k^* + \mathbf{x}^\top \boldsymbol{\gamma}_k^* + \Gamma \|\boldsymbol{\gamma}_k^*\|_*.$$

*Because $(\rho_k^*, \boldsymbol{\gamma}_k^*)$ is the minimum for the fixed $\boldsymbol{\pi}^{(k)}$, we have*

$$v_k(\mathbf{z},\mathbf{x}) = \min_{\rho,\boldsymbol{\gamma}} \left\{ \rho + \mathbf{x}^\top \boldsymbol{\gamma} + \Gamma \|\boldsymbol{\gamma}\|_* \mid \rho + \mathbf{x}_s^\top \boldsymbol{\gamma} \geq (\boldsymbol{\pi}^{(k)})^\top (\mathbf{h}_s - \mathbf{T}_s \mathbf{z}), \ \forall s \right\}.$$

*On the other hand, in formulation (9) $\mathcal{Q}(\mathbf{z},\mathbf{x})$ is defined as the maximum over all $\boldsymbol{\pi} \in \Pi$ of precisely the same inner problem. Hence*

$$\mathcal{Q}(\mathbf{z},\mathbf{x}) \geq v_k(\mathbf{z},\mathbf{x}). \tag{$\star$}$$

*Inequality ($\star$) together with the definition of $v_k(\mathbf{z},\mathbf{x})$ yields*

$$\alpha \geq \mathcal{Q}(\mathbf{z},\mathbf{x}) \implies \alpha \geq \rho_k^* + \mathbf{x}^\top \boldsymbol{\gamma}_k^* + \Gamma \|\boldsymbol{\gamma}_k^*\|_*,$$

*which is exactly (17a). Because $k$ was arbitrary, all cuts in (17) are valid for every $(\mathbf{z},\alpha)$ that satisfies (18b). Therefore the master problem that contains the stored columns remains a relaxation of the original upper-level problem irrespective of the new covariate vector $\mathbf{x}$.*

**Remark:** The adaptive nature of the Contextual CCG algorithm underscores its paramount significance, as it has the potential to alleviate computational burden in situations where optimization needs sequential updates of new contextual information, such as in online optimization applications. Additionally, it proves beneficial in circumstances where alternative algorithms result in prohibitively long computational times, striking a balance between the computational load to solve the master and the number of iterations required for convergence.

*Finally, the Contextual CCG algorithm must be executed by solving the master problem (16) instead of (14) at line 4, using a set $\{\boldsymbol{\pi}^{(k)}\}_{k \in [K]}$ that defines the pre-computed cuts, and by replacing (14) with (16) at line 8.*

# 5 Extensions

*In this section, we present extensions for the recourse objective uncertainty case and the more general case involving both right-hand-side and objective uncertainty.*

## 5.1 Recourse problem with objective uncertainty

*In this section, we present a tractable reformulation of the problem (3) by assuming that uncertainty is only present in the recourse vector cost ($\tilde{\mathbf{y}} := \mathbf{q}$). The primal formulation of the recourse problem is given by:*

$$Q(\mathbf{z}, \tilde{\mathbf{y}}) = \min_{\mathbf{u} \geq 0} \tilde{\mathbf{y}}^\top \mathbf{u} \tag{19a}$$

$$\text{s.t.} \quad \mathbf{W}\mathbf{u} = \mathbf{h} - \mathbf{T}\mathbf{z}. \tag{19b}$$

**Proposition 7** *If uncertainty affects only the objective function of the lower-level problem in (4), then problem (3) is equivalent to*

$$\min_{\mathbf{z}, \mathbf{u}, \rho, \gamma} \quad \mathbf{c}^\top \mathbf{z} + \rho + \mathbf{x}^\top \gamma + \Gamma \|\gamma\|_*$$
$$\text{s.t.} \quad \mathbf{W}\mathbf{u} = \mathbf{h} - \mathbf{T}\mathbf{z},$$
$$\rho + \mathbf{x}_s^\top \gamma \geq \mathbf{q}_s^\top \mathbf{u}, \qquad \forall s \in [S], \tag{20}$$
$$\mathbf{u} \geq 0,$$

*where $\|\cdot\|_*$ denotes the dual norm.*

**Proof 8** *Proof of Proposition 7 If $Y_\Gamma(\mathbf{x}) \ni \tilde{\mathbf{y}} = \sum_{s \in [S]} \theta_s \mathbf{q}_s$, then the recourse problem is:*

$$Q(\mathbf{z}, \tilde{\mathbf{y}}) = \min_{\mathbf{u} \geq 0} \sum_{s \in [S]} \theta_s \left( \mathbf{q}_s^\top \mathbf{u} \right) \tag{21a}$$

$$\text{s.t.} \quad \mathbf{W}\mathbf{u} = \mathbf{h} - \mathbf{T}\mathbf{z}. \tag{21b}$$

*Therefore, we have:*

$$Q(\mathbf{z}, \mathbf{x}) = \max_{\varepsilon, \theta \geq 0} \left\{ \min_{\mathbf{u} \geq 0} \left\{ \sum_{s \in [S]} \theta_s \left( \mathbf{q}_s^\top \mathbf{u} \right) \,\middle|\, \mathbf{W}\mathbf{u} = \mathbf{h} - \mathbf{T}\mathbf{z} \right\} \,\middle|\, \begin{array}{l} \sum_{s \in [S]} \theta_s = 1 \\ \sum_{s \in [S]} \theta_s \mathbf{x}_s - \varepsilon = \mathbf{x} \\ \|\varepsilon\| \leq \Gamma \end{array} \right\}, \tag{22}$$

*where $\theta = (\theta_1, \ldots, \theta_S)^\top$. From the dual formulation of the recourse function, which is the inner problem in (22), we obtain:*

$$Q(\mathbf{z}, \mathbf{x}) = \max_{\varepsilon, \theta \geq 0} \left\{ \max_{\pi} \left\{ \pi^\top \left( \mathbf{h} - \mathbf{T}\mathbf{z} \right) \,\middle|\, \mathbf{W}^\top \pi \leq \sum_{s \in [S]} \theta_s \mathbf{q}_s \right\} \,\middle|\, \begin{array}{l} \sum_{s \in [S]} \theta_s = 1 \\ \sum_{s \in [S]} \theta_s \mathbf{x}_s - \varepsilon = \mathbf{x} \\ \|\varepsilon\| \leq \Gamma \end{array} \right\}. \tag{23}$$

*Formulation (23) can be reformulated as the following single-level maximization problem:*

$$\max_{\pi, \varepsilon, \theta \geq 0} \quad \pi^\top \left( \mathbf{h} - \mathbf{T}\mathbf{z} \right)$$
$$\text{s.t.} \quad \mathbf{W}^\top \pi \leq \sum_{s \in [S]} \theta_s \mathbf{q}_s, \quad : \mathbf{u}$$
$$\sum_{s \in [S]} \theta_s = 1, \quad : \rho \tag{24}$$
$$\sum_{s \in [S]} \theta_s \mathbf{x}_s - \varepsilon = \mathbf{x}, \quad : \gamma$$
$$\|\varepsilon\| \leq \Gamma.$$

*From the dual of problem (24), we obtain the following dual problem:*

$$\min_{\mathbf{u}, \rho, \gamma} \quad \rho + \mathbf{x}^\top \gamma + \Gamma \|\gamma\|_*$$
$$\text{s.t.} \quad \mathbf{W}\mathbf{u} = \mathbf{h} - \mathbf{T}\mathbf{z},$$
$$\rho + \mathbf{x}_s^\top \gamma \geq \mathbf{q}_s^\top \mathbf{u}, \quad \forall s \in [S], \tag{25}$$
$$\mathbf{u} \geq 0.$$

*By combining this dual problem with the first-stage problem, we obtain the following overall single-level opti-*

*mization problem:*

$$\min_{\mathbf{z},\mathbf{u},\rho,\gamma} \quad \mathbf{c}^\top \mathbf{z} + \rho + \mathbf{x}^\top \gamma + \Gamma \|\gamma\|_*$$

$$\text{s.t.} \quad \mathbf{W}\mathbf{u} = \mathbf{h} - \mathbf{T}\mathbf{z},$$
$$\rho + \mathbf{x}_s^\top \gamma \geq \mathbf{q}_s^\top \mathbf{u}, \qquad \forall s \in [S], \quad (26)$$
$$\mathbf{u} \geq 0,$$
$$\mathbf{z} \in Z.$$

## 5.2 Model with full data uncertainty

*In the presence of joint uncertainty in the objective coefficients, right-hand side, and recourse matrix, i.e., $\mathbf{y} := (\mathbf{q}, \mathbf{h}, \mathbf{T}, \mathbf{W})$, the worst-case function $\mathcal{Q}(\mathbf{z}, \mathbf{x})$ is defined by the following bilevel problem:*

$$\max_{\varepsilon, \theta \geq 0} \min_{\mathbf{u} \geq 0} \left\{ \sum_{s \in [S]} \theta_s \left( \mathbf{q}_s^\top \mathbf{u} \right) \,\middle|\, \mathbf{W}\mathbf{u} = \sum_{s \in [S]} \theta_s \left( \mathbf{h}_s - \mathbf{T}_s \mathbf{z} \right) \right\} \quad (27\text{a})$$

$$\text{s.t.} \quad \sum_{s \in [S]} \theta_s = 1 \quad (27\text{b})$$

$$\sum_{s \in [S]} \theta_s \mathbf{x}_s - \varepsilon = \mathbf{x} \quad (27\text{c})$$

$$\|\varepsilon\| \leq \Gamma. \quad (27\text{d})$$

*In this case, problem (3) can be addressed using the CCG method. The oracle problem is defined by (27), and constraint (15b) is replaced by $\alpha \geq \left(\mathbf{q}^{(k)}\right)^\top \mathbf{u}_k$, $\forall k \in [K]$, where $\mathbf{q}^{(k)} = \sum_{s \in [S]} \theta_s^{(k)} \mathbf{q}_s$. Here, $\theta^{(k)}$ is the solution of (27) corresponding to the suboptimal solution $\mathbf{z}^{(k)}$ at each iteration $k \in [K]$. However, the bilevel problem (27) can still be challenging to solve, even for moderately sized recourse problems.*

# 6 Application in the online energy and reserve scheduling problem

*In this section, we demonstrate the effectiveness of our framework through a widely studied application in online energy operations: the hour-ahead economic dispatch problem for energy and reserve scheduling. This problem plays a critical role in the operation of electricity markets, where increasing penetration of intermittent renewable energy sources introduces significant operational uncertainty. As the share of non-dispatchable generation grows, system operators require more sophisticated predictive and prescriptive models to manage uncertainty and maintain grid reliability.*

*Robust Optimization offers key advantages in this setting, providing transparent, reproducible, and compliance-friendly decision support tools. In contrast to stochastic programming approaches, which may produce sample-dependent solutions, robust methods naturally yield decisions with strong out-of-sample performance guarantees—an important consideration for regulatory and operational acceptance in power system applications.*

*Motivated by industry practices, we develop a case study based on the hour-ahead economic dispatch problem under uncertainty in non-dispatchable renewable energy generation (REG). We formulate it as a two-stage conditional-robust optimization problem: in the first stage, energy and reserve schedules are determined; in the second stage, a re-dispatch plan is computed for the worst-case realization within the DDCU set, subject to the first-stage decisions. We begin by analyzing the properties of the DDCU set for REG, followed by a presentation of our approach's numerical results and performance assessment.*

*Finally, we emphasize that the energy application problem has uncertainty only in the vector $\mathbf{h}$, i.e., $\mathbf{y} := \mathbf{h}$. Therefore, we conducted the computational experiments using the right-hand-side solution methodology described in Section 4, based on the classical P-Bilevel formulation (5) and the $D_h$-Bilevel formulation (11) of the oracle problem.*

## 6.1 Conditional Uncertainty Set Driven by Renewable Energy Generation (REG) Data

*In generation scheduling with a high penetration of renewable-based generation, uncertain data include the power output for each renewable generator and time period within the scheduling horizon. In the current state-of-the-art, the characterization of this uncertainty affects a broad class of predictive models that utilize temporal cross-lagged dependencies and categorical data (e.g., month of the year, day of the week, hour of the day, weather forecasts) as features.*

*In this work, we employ historical hourly profiles and features like cross-lagged dependencies and the hour-of-the-day indicator as scenarios. This allows us to endogenously consider the multidimensional dependencies in true hourly REG and bypass the predictive step when quantifying uncertainty. Mathematically, a scenario $s \in [S]$ consists of the coordinates $(\mathbf{x}_s, \mathbf{y}_s)$, where $\mathbf{x}_s$ represents the observed cross-lagged power output of a renewable generator for each period or specific hour of the day, and $\mathbf{y}_s$ represents a historical hourly profile of REG. The capability of the DDCU set to capture temporal dependencies is illustrated in Figure 4, where two-dimensional projections of the multidimensional DDCU set are shown in yellow for the data of the case study examined in Section 6.2.*

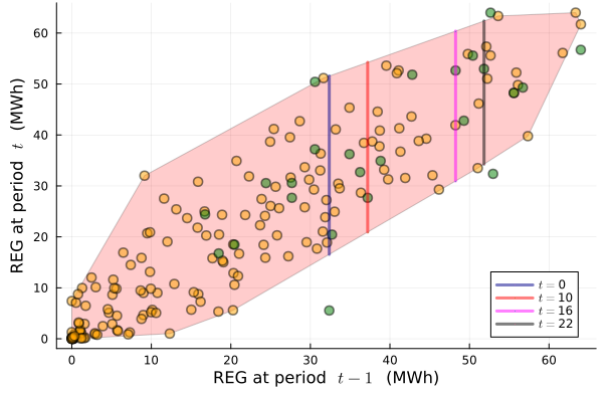

Figure 4: Example of temporal cross-lagged dependences.

*On the x-axis of Figure 4, we show the power output lagged by $t-1$, and on the y-axis, the power output of the renewable generator at some bus of the system for period t, covering all periods $t = 1, \ldots, 24$. The green bullet points represent the two-dimensional projections of the power output observed at period t, explained by the power output at $t-1$, which can fall outside the convex hull formed by the observations. The lines represent two-dimensional projections of the DDCU set corresponding to observations of the contextual information for time periods $t = 1, 10, 16, 22$. We can observe how the DDCU set changes as contextual information is revealed at each period of system operation.*

*As discussed in Section 2, the way the DDCU set is constructed makes it adaptive to new observations of contextual information. In particular, when the DDCU set is built using historical hourly REG profile data, this adaptiveness gives rise to a rolling-window scheme that reflects how power systems are typically operated in practice by the electrical industry. This approach will be used in the numerical experiments reported in the following sections.*

*Finally, to run our computational experiments, we calibrate the parameter $\Gamma$ as an $\delta$-percent increase over the $\Gamma_0$ value, i.e., $(1 + \delta)\Gamma_0$. We set $\varepsilon = 10\%$, that is, $\Gamma$ is 10% larger than $\Gamma_0$ for each new observation of the contextual information at each time period t. As each new observation of the contextual information is revealed at time t, we compute its distance to the polyhedron defined by the convex hull of the historical data within the rolling window, which is used to define $\Gamma_0$.*

## 6.2 Numerical results

We apply our methodology to three distinct power system test cases: the IEEE 118-bus test system, a 600-bus amalgamation formed by combining two instances of the IEEE 300-bus system, and the 10,000-bus Grid Optimization Competition (GOC) test case [Babaeinejadsarookolaee et al., 2021]. To assess the efficacy of the proposed framework, we investigate wind-related uncertainties, using data from the Global Energy Competition (GEFCom) dataset [Hong et al., 2014, Hong et al., 2016]. The wind generation data underwent rescaling from its original configuration involving 10 buses and 10 wind farms. Our simulations were executed using Gurobi 10.0.0 within the JuMP (Julia 1.8) modeling language, running on an Xeon E5-2680 processor operating at 2.5 GHz with 128 GB of RAM. The optimality gap of Gurobi was set at 0.01%

We benchmark our framework against the Data-Driven Robust Optimization model proposed in [Velloso et al., 2020], hereafter referred to as the unconditional case, for the energy and reserve scheduling problem in an hour-ahead economic dispatch instance. First, we discuss the economic benefits of our proposed methodology in comparison to the unconditional case. Then, we analyze the computational performance of the algorithmic schemes described in Section 4 for the proposed application.

### 6.2.1 Economic and reliability analysis

We evaluated the energy application problem as defined in (28), whose detailed problem formulation, can be found in Appendix A. This evaluation served as a benchmark for our framework against the unconditional case, incorporating an out-of-sample economic and reliability analysis.

To initiate this process, we created the DDCU set using a variety of covariate data from different training samples. We analyzed 840 hourly renewable energy generation (REG) profiles, reflecting the power output of 10 renewable generators, to form our uncertainty scenarios. For covariate scenarios, we incorporated elements typically used in auto-regressive (AR) predictive models for REG time series. This included the lagged power outputs at $t-1$ ($AR_1$) and the hour of the day, which was represented as a dummy variable. These datasets were specifically developed for the IEEE 118-bus test system.

Following the creation of these datasets, we proceeded to solve the energy and reserve scheduling problem in a sequential manner for each period $t$, within a time horizon extending from $t = 0$ to $t = 24$. This comprehensive approach allowed for a thorough assessment of the performance and reliability of our framework under varying and realistic operational conditions.

It is worth mentioning that our DDCRO framework exhibits adaptability in two significant ways. First, as the problem exhibits right-hand side uncertainty only, our D-Bilevel reformulation of the oracle (9)—in addition to ensuring tractability—integrates explicit contextual information into the master for the Contextual CCG algorithm. This integration not only captures the complete structure of the DDCU set but also facilitates the updating of the predictive-prescriptive model (3) with new covariate observations. Crucially, this update occurs without the need for model re-training, allowing for the commencement of computations with an effective sub-optimal solution.

Second, the adaptability of our framework is further enhanced by the method of constructing the DDCU set. This construction allows the predictive-prescriptive model to be dynamically updated, utilizing a rolling-window scheme, one step ahead of the last S historical observations of both uncertainty and contextual information. This approach, as outlined in sources such as [Fernandes et al., 2016, Velloso et al., 2020], ensures that the model remains responsive and current, effectively incorporating the latest available data for more accurate and timely decision making.

For the energy application problem, we performed an out-of-sample analysis using an adaptive one-step-ahead forecast for the power output of renewable generators within a rolling window of 720 periods (one month). With the DDCU set in hand, we sequentially solved instances of the energy and reserve scheduling problem for each period of the day-ahead. The total cost for that day was obtained from the solutions of these instances. We then calculated the infeasibility of the generation dispatch linked to the produced generation schedules and the actual REG scenario observed for the day ahead to determine the robustness of the solution. To achieve this, we solved a single-period version of the inner minimization problem in (28a)-(28q) for each *hour of the day, following standard industrial practice. Consequently, we developed out-of-sample statistics for robustness based on the hourly amounts of load shedding and REG leakage that were generated.*

| Model | Avg. total cost ($10^4$\$) | Diff. avg. total cost (%) | LOLP (%) | PWS (%) |
|---|---|---|---|---|
| Unconditional | 134,0 | - | 0,14 | 3,61 |
| Conditional dummy | 132,2 | 1,3 | 0,14 | 7,64 |
| Conditional $AR_1$ | 125,7 | 6,2 | 0,83 | 7,78 |
| Conditional $AR_1$ + dummy | 124,7 | 6,9 | 0,69 | 10,0 |

Table 1: Results for the 118-bus system.

*Column 1 of Table 1 provides the average total cost of day-ahead operations for the IEEE 118-bus test system. Column 2 displays the percentage difference between the unconditional and conditional cases. Column 3 presents a reliability index referred to as the Loss of Load Probability (LOLP), which is defined as the fraction of hours with load shedding exceeding 0.1% of the system load. Similarly, column 4 reports the Probability of Wind Spillage (PWS), defined as the fraction of hours in which wind spillage exceeds 0.1% of the available wind power generation.*

*The first observation is that, in the unconditional case, the precision of the DDCU set (defined by the training sample size) allows moving in just one direction: towards high reliability (low risk reflected in small LOLP values) with a significantly high cost. In the conditional case, the precision of the DDCU set steers the solution towards lower reliability with lower cost. Moreover, we can adjust the precision level of the DDCU set by increasing the δ-percent of the $Γ_0$ value to calibrate Γ, aiming to achieve a solution with significantly lower cost and LOLP values still below 1%, representing a low level of risk.*

*We observe that the PWS is significantly greater than the LOLP, especially in the case of the contextual robust optimization model. This is because the contextual robust optimization model is designed in a more informed and expert manner compared to its unconditional counterpart, particularly when out-of-sample scenarios are revealed. The contextual robust model dispatches a base set of thermal generators that can meet system demand without incurring the highest operational costs, thus avoiding load shedding in cases of low power output of renewable generators. In scenarios with high power output from renewable generators, it is possible to re-dispatch a few generators and allow for some renewable spillage, while still maintaining a low overall system cost. In contrast, due to its inherent conservatism, the unconditional robust optimization model commits to the highest-cost dispatch strategy to avoid load shedding, as it prioritizes scenarios with low renewable power output.*

### 6.2.2 Computational performance analysis

*In this section, we first emphasize the critical need to adopt the D-Bilevel formulation (9) over the P-Bilevel approach by examining the computational load of the classical CCG algorithm under both oracles. We then benchmark our proposed Contextual CCG algorithm, presented in Section 4, against the classical CCG with the P-Bilevel oracle to assess their computational efficiency. The analysis focuses on their application to sequentially solving the energy dispatch problem (28) for scheduling energy and reserves over periods $t = 1, …, 24$..*

*For this purpose, we utilized a training data set comprising 840 hourly REG profiles, along with lagged $t − 1$ power outputs from 10 renewable generators (serving as covariates), to construct the DDCU set for the energy application problem as detailed in Section 6.1. We conducted a computational experiment by applying the energy problem to the referenced test systems, employing both the classical CCG algorithm and its adaptations for D-Bilevel and P-Bilevel formulations of the oracle. For this purpose, we utilized the BilevelJuMP.jl toolbox, integrated into the Julia programming language, as documented in [Garcia et al., 2023]. Both bilevel optimization problems were approached using the big-M method ([Fortuny-Amat and McCarl, 1981]) and incorporated within the mode options of the toolbox.*

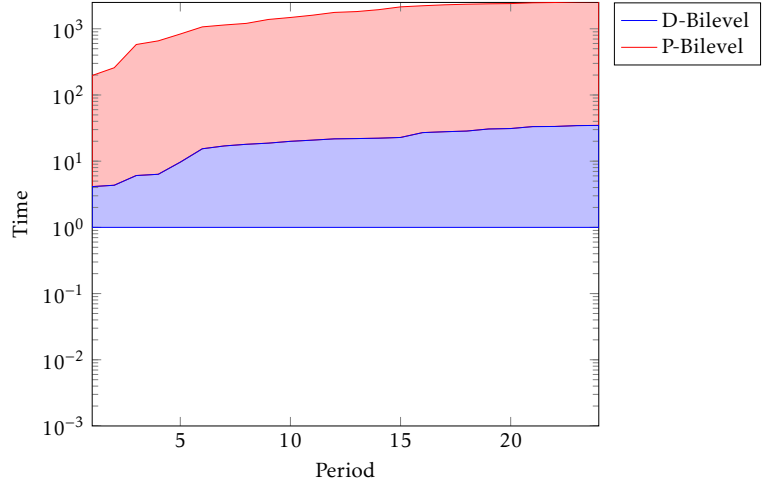

Figure 5: Cumulative convergence duration (logarithmic scale) of classical CCG methods, measured in minutes, over multiple periods for the 600-bus system.

*Specifically, for the IEEE 118-bus test system, the algorithm employing the P-Bilevel formulation of the oracle demonstrated superior performance compared to the $D_h$-Bilevel formulation. This was primarily because the lower-level complementary slackness condition in the P-Bilevel formulation was more manageable due to the lower complexity of the IEEE 118-bus system. Conversely, in the case of the more complex 10,000-bus test system, the P-Bilevel approach resulted in prohibitively long computational times. Here, the D-Bilevel formulation proved to be crucial, rendering the bilevel optimization problem tractable. We further investigated this phenomenon using an intermediate complexity instance, such as the 600-bus test system, to understand the scalability and applicability of these formulations.*

*As shown in Figure 5, the cumulative convergence duration over the 24 time periods for the algorithm using the D-Bilevel formulation is well below $10^2$, whereas the P-Bilevel formulation exceeds $10^3$ by a significant margin. This stark difference underscores the necessity of developing the $D_h$-Bilevel formulation. When confronted with increasing complexity and scalability challenges in the lower-level problem of the P-Bilevel formulation, resorting to the $D_h$-Bilevel approach becomes essential to maintain the tractability of the oracle.*

*Figure 6 compares the convergence time of the adaptive Contextual CCG algorithm—which leverages the converged master from the previous operation period to sequentially solve the next problem—against the classical CCG algorithm that employs the standard P-Bilevel oracle formulation. It is worth noting that in this classical formulation, which represents the state of the art for right-hand-side uncertainty problems, the contextual information is endogenously embedded, thereby preventing the adaptive reuse of the converged master from previous operation periods.*

*This result demonstrates that our $D_h$-Bilevel formulation of the oracle problem not only renders the problem tractable but also enables the development of a novel adaptive CCG algorithm with explicit contextual information, which significantly outperforms the benchmark for right-hand-side uncertainty in the online energy application problem.*

*Table 2 presents the average computational time required for the convergence of both algorithms, averaged across all time periods $t = 1, \ldots, 24$. For each time period, we calculate the cumulative time required to solve both the master problem and the oracle for incorporating a new dual vertex $\pi^{(k)}$ or primal vertex $\theta^{(k)}$ into the master problem during each iteration k. The data for each test system is organized into three columns: the first column details the average time to solve the master problem, the second column indicates the average time for solving the oracle, and the third column combines these to show the total average time required.*

*As previously mentioned, the P-Bilevel formulation incurs prohibitively long computational times. The Contextual CCG algorithm achieved the best results. The explicit integration of contextual information renders this algorithm adaptive, enabling the master problem to be initialized with a set $\{\pi^{(k)}\}_{k \in [K]}$ of pre-computed*

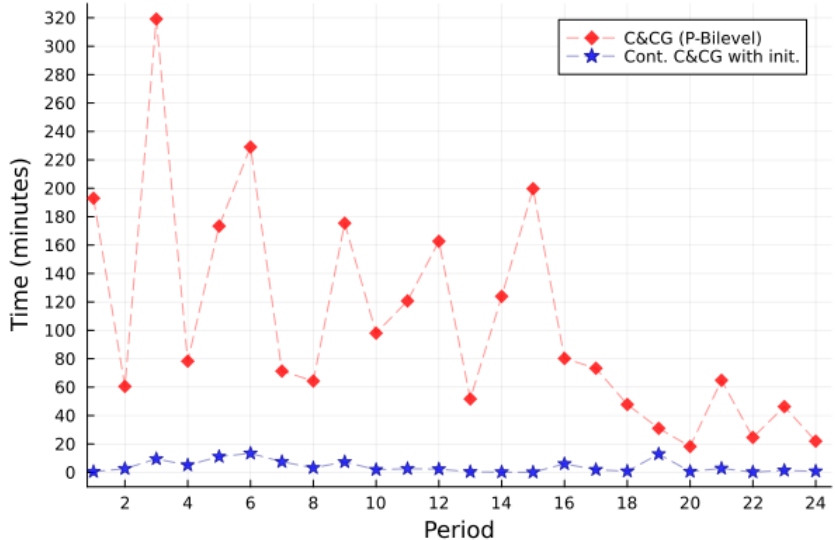

Figure 6: Computation time: classical CCG (P-Bilevel oracle) vs. Contextual CCG with period-wise master initialization.

| Method | Oracle formulation | 600-bus test syst. | | | 10,000-bus test syst. | | |
|---|---|---|---|---|---|---|---|
| | | Avg. time master ($10^{-3}$) | Avg. time oracle | Avg. total time | Avg. time master | Avg. time oracle | Avg. total time |
| CCG | P-Bilevel | 4,7 | 105,4 | 105,4 | - | - | - |
| Contextual CCG | $D_h$-Bilevel | 12,0 | 4,0 | 4,0 🏆 | 9,3 | 24,3 | 33,6 🏆 |

Table 2: Average computational time (in minutes) required for the convergence of each algorithm in solving the energy application problem for the 600-bus and 10,000-bus test systems.

vertices for any new observation of contextual information. This leads to a more effective starting point for the algorithm, providing a sub-optimal—but of good quality—solution that simplifies the feasibility adjustment process.

To analyze the master–oracle relationship of the proposed Contextual CCG algorithm, we focus on the more complex 10,000-bus test system. Figure 7 shows a stacked bar chart reporting the time required to solve both the master and the oracle across 24 periods of a standard operating day. The hour-ahead energy and reserve scheduling problem is inherently challenging due to the operational dynamics of power systems, with some periods proving more difficult than others. Typically, the computational load of the bilevel oracle problem, as defined by the $D_h$-Bilevel formulation, exceeds that of the master problem, and this discrepancy becomes especially pronounced during challenging periods, such as at $t = 12$. The results also reveal a characteristic trade-off between the master and oracle problems. While the master may find a near-optimal solution—usually incurring significant computational load to exploit precomputed cuts from the contextual information—the oracle can often certify the optimality of such suboptimal solutions more quickly. Conversely, when the master yields poor suboptimal solutions, which typically require little computation time, the oracle must perform more iterations or may even struggle to certify optimality, ultimately leading to a significant computational burden.

The Contextual CCG algorithm, thanks to its adaptive nature and a well-chosen initial sub-optimal solution, significantly eases the complexity of the bilevel optimization problem defined by the $D_h$-Bilevel formulation. This advantage is particularly evident in challenging periods, underscoring the efficiency of the algorithm and adaptability to manage complex operational circumstances in power systems.

It is worth emphasizing that the Contextual CCG method uses a lighter master problem, which reduces the computational load. Moreover, it incorporates pre-computed cuts that are valid for any realization of contextual information. Empirically, this can be explained by the fact that similar realizations of contextual

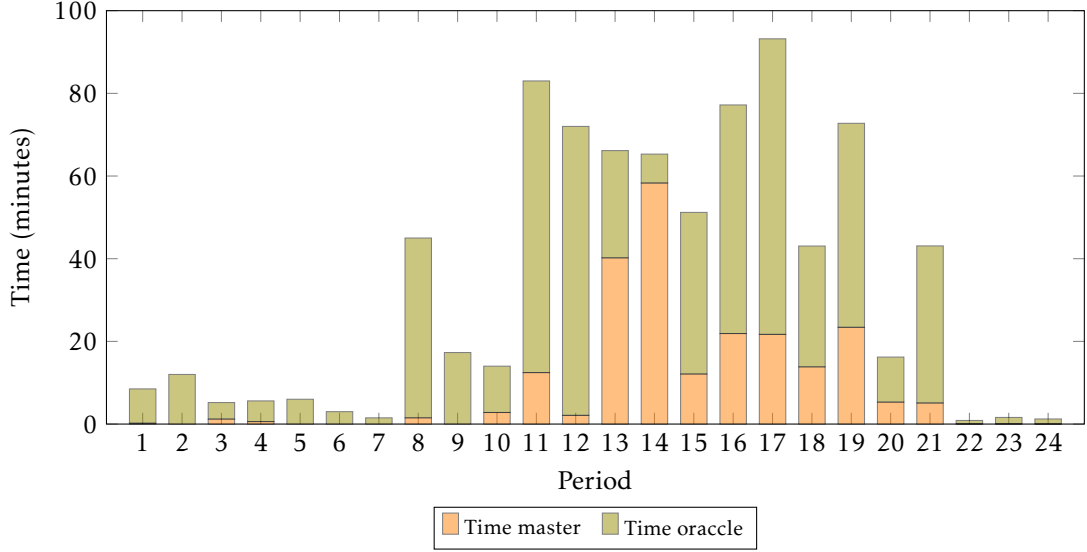

Figure 7: Convergence time of the proposed Contextual CCG algorithm on the 10,000-bus system.

*information should lead to master problems with similar sets of cuts. Thus, when a new observation of contextual information closely resembles a previous one—whose cuts have already been incorporated—the Contextual CCG algorithm may require only a few iterations, or even just one, to verify or assert optimality by querying the oracle.*

*In summary, the Contextual CCG method inherits key advantages from its counterparts: it features a lightweight master problem can potentially assert optimality in few iterations, as in the classical CCG algorithm, allowing us to solve large-scale problems.*

# 7 Conclusions

*We introduce a novel Data-Driven Contextual Robust Optimization (DDCRO) framework that fundamentally advances the integration of online contextual information into data-driven robust optimization models. Existing end-to-end contextual robust formulations are typically trained for specific contexts and require full retraining whenever new contextual information becomes available. Such retraining hinders their applicability to online decision-making problems, especially those requiring frequent re-optimization, such as multiperiod energy scheduling. The proposed DDCRO framework overcomes this limitation by jointly embedding contextual learning and prescriptive decision-making within a unified optimization architecture.*

*At the core of the framework lies the Data-Driven Contextual Uncertainty (DDCU) set, a polyhedral uncertainty set constructed directly from empirical data and conditioned on contextual covariates. This DDCU can be interpreted as a data-driven nonparametric regression model embedded within the optimization problem, where contextual information serves as regressors that locally adapt the geometry of the uncertainty set. In this sense, the DDCU set plays the dual role of a predictive model and a robust uncertainty representation. It captures how uncertainty evolves with context without imposing parametric assumptions, while ensuring that the resulting optimization problem remains tractable. This nonparametric conditioning process allows the optimization model to "learn" directly from data in an end-to-end fashion, bypassing the need for a separate predictive step.*

*From a methodological perspective, our framework departs from traditional robust optimization paradigms in several important ways. First, it eliminates the need for retraining by embedding contextual information directly into the cut generation process. This is due to a novel reformulation of the second and third levels of the proposed DDCRO problem, which yields a Contextually-Adapted Column-and-Constraint Generation algorithm suitable for online decision-making processes. Thus, because the contextual variables are explicitly*

*considered in the cut expressions, cuts derived from previous periods, to prior contexts, remain valid for new information. The result is a decomposition algorithm that not only accelerates convergence but also adapts dynamically to evolving data streams. This property is particularly advantageous in energy scheduling, where decision problems are repeatedly solved every hour, or few minutes, under continuously changing uncertainty conditions.*

*Our computational experiments in the context of energy and reserve scheduling with uncertain renewable generation offer compelling evidence of these advances. The test systems, including the IEEE 118-bus, a 600-bus composite of two IEEE 300-bus systems, and the 10,000-bus Grid Optimization Competition case provide compelling evidence of scalability and gains in comparison to the state-of-the-art benchmark used to solve the proposed class of problem. Across all cases, the DDCRO framework consistently outperforms the classical data-driven robust optimization benchmark, both in economic performance and computational efficiency.*

*In economic terms, the conditional (contextual) models achieved up to a 7% reduction in total operating costs relative to unconditional robust formulations. These gains arise because contextual conditioning reduces conservatism: the uncertainty set adapts to the prevailing covariate information (such as time of day, forecasted conditions, or lagged renewable generation), thereby excluding unrealistic or irrelevant worst-case realizations. Consequently, the resulting dispatch schedules are both cost-effective and reliable, maintaining loss-of-load probabilities under 1% while preserving operational feasibility across out-of-sample scenarios. This outcome reflects a shift from overly conservative decisions to contextually aware data-driven robustness, an essential balance for modern power systems operating with high shares of intermittent generation.*

*From a computational perspective, the proposed Contextual CCG algorithm provides substantial efficiency gains. For large-scale systems, the algorithm leveraging the D-Bilevel oracle formulation achieves convergence times one to two orders of magnitude faster than the traditional P-Bilevel implementation. These results confirm that the reformulation approach used in this work and the resulting reuse of valid cuts directly contribute to the computational tractability of large-scale practical problems. Notably, in the 10,000-bus case, the classical formulation and CCG approach become impractically slow, while the D-Bilevel-based Contextual CCG retains tractability.*

*Beyond the specific application in power system operations, the DDCRO framework contributes a general methodology for contextually adaptive optimization under uncertainty. The ability to construct uncertainty sets conditioned on observed covariates opens new avenues for applying robust optimization in domains where environmental or market conditions evolve dynamically. Potential applications include financial portfolio management, logistics planning under demand uncertainty, healthcare resource allocation, and industrial production scheduling. In all these cases (ranging from macroeconomic indicators to weather forecasts or production metrics), the nonparametric regression based on the DDCU set plays a decisive role in shaping uncertainty distributions for decision-making under uncertainty. The DDCRO methodology provides a principled way to incorporate these variables without sacrificing tractability or interpretability.*

*From a policy and operational standpoint, our findings have clear implications for the energy sector. As power systems transition toward higher shares of renewable generation, operators must adopt decision-support tools capable of balancing reliability, economic efficiency, and computational feasibility. However, accountability, transparency, and reproducibility are often overlooked in academic publications, yet crucial and determinant in the selection of a practical decision tool (see [Dias Garcia et al., 2025] for an in-depth discussion). In this sense, the proposed DDCRO framework constitutes a robust framework; therefore, it is a deterministic method (not sample-dependent as classical stochastic approaches based on SSA) that is reproducible and based on real data, characteristics of great value for system operators to justify their choice.*

*The present work lays the groundwork for several promising research directions. Extending DDCRO to multi-stage and dynamic settings, where uncertainty unfolds sequentially and decisions must adapt over longer horizons, constitutes a natural step. Further exploration of distributionally robust variants of the DDCU set, potentially in combination with the multistage framework, is also of interest. Finally, applying the proposed method to other industrial contexts would further demonstrate its versatility and practical relevance.*

# A    Online energy application problem

*In this section, we present the online energy application problem, which serves as a benchmark for evaluating the effectiveness of our proposed framework. Our methodology starts with the establishment of a clear and comprehensive set of notations, specifically tailored to this problem. Once these notations are defined, we will meticulously detail the mathematical formulation of the problem. The formulation encompasses various key components, including objective functions, constraints, and decision variables, each playing an essential role in the overall model.*

## A.1    Notation

*In this section, we present a notation list for the online energy application problem. Bold symbols are allocated to matrices (uppercase) and vectors (lowercase).*

*Constants*

| | |
|---|---|
| $\mathbf{G}$ | Thermal generator-bus incidence matrix. |
| $\mathbf{A}$ | Line-bus incidence matrix. |
| $\mathbf{D}_t$ | Vector of nodal consumption at period t. |
| $\mathbf{B}$ | Angle-to-flow matrix. |
| $\mathbf{E}$ | Renewable generator-bus incidence matrix. |
| $\overline{\mathbf{P}}$, $\underline{\mathbf{P}}$ | Vector of the maximum and minimum generation limits of thermal units. |
| $\mathbf{R}^{dn}$, $\mathbf{R}^{up}$ | Vector of downward and upward limits for corrective actions of thermal units within each period. |
| $\overline{\mathbf{F}}$ | Maximum power flow. |
| $\overline{\mathbf{y}}$ | Vector of the expected REG. |
| $\mathbf{c}$ | Vector of fuel costs of thermal generators. |
| $\mathbf{c}^{dn}$, $\mathbf{c}^{up}$ | Vectors of cost rates for down- and up-spinning reserves. |

*Variables*

| | |
|---|---|
| $\mathbf{p}_t$ | Vector of nominal power generation at time period t. |
| $\mathbf{r}_t^{dn}$, $\mathbf{r}_t^{up}$ | Down- and up-spinning reserves at time period t. |
| $\beta_t$ | Phase angles at time period t. |
| $\mathbf{f}_t$ | Line power flows at time period t. |
| $\Delta_t$ | Decision vector representing the re-dispatch of energy generation at time period t. |
| $\beta_t^{wc}$ | Decision vector representing the conditioned worst-case nodal phase angles at time period t. |
| $\mathbf{f}_t^{wc}$ | Decision vector representing the conditioned worst-case line power flows at time period t. |
| $\mathbf{s}_t^+$, $\mathbf{s}_t^-$ | Decision vectors representing the REG spillage and load shedding at time period t. |

## A.2    Mathematical formulation of the online energy application problem

*Using previous work on robust energy and reserve scheduling [Moreira et al., 2014, Cobos et al., 2018, G. Cobos et al., 2018] and well-known industry practices [Arroyo and Galiana, 2005, Galiana et al., 2005], we consider a base case dispatch for energy and reserve that ensures $\alpha$-feasibility for all scenarios within the DDCU set. The mathematical formulation for the co-optimized hour-ahead energy and reserve problem is presented below:*

$$\min_{\substack{\mathbf{p}_t, \mathbf{r}_t^{up}, \mathbf{r}_t^{dn}, \\ \mathbf{f}_t, \beta_t, \alpha}} \mathbf{c}^\top \mathbf{p}_t + (\mathbf{c}^{up})^\top \mathbf{r}_t^{up} + (\mathbf{c}^{dn})^\top \mathbf{r}_t^{dn} + \alpha \tag{28a}$$

$$\text{s.t.} \quad \mathbf{G}\mathbf{p}_t + \mathbf{A}\mathbf{f}_t = \mathbf{D}_t - \mathbf{E}\overline{\mathbf{y}} \tag{28b}$$

$$\mathbf{f}_t = \mathbf{B}\beta_t, \tag{28c}$$

$$\underline{\mathbf{P}} + \mathbf{r}_t^{dn} \leq \mathbf{p}_t \leq \overline{\mathbf{P}} - \mathbf{r}_t^{up}, \tag{28d}$$

$$-\overline{\mathbf{F}} \leq \mathbf{f}_t \leq \overline{\mathbf{F}}, \tag{28e}$$

$$0 \leq \mathbf{r}_t^{up} \leq \mathbf{R}^{up}, \tag{28f}$$

$$0 \leq \mathbf{r}_t^{dn} \leq \mathbf{R}^{dn}, \tag{28g}$$

$$\mathbf{p}_t + \mathbf{r}_t^{up} \leq \mathbf{R}^{up} + \mathbf{p}_{t-1} \tag{28h}$$

$$\mathbf{p}_t - \mathbf{r}_t^{dn} \geq \mathbf{p}_{t-1} - \mathbf{R}^{dn} \tag{28i}$$

$$\mathcal{Q}(\mathbf{z}, \mathbf{x}_t) \leq \alpha \tag{28j}$$

$$\mathcal{Q}(\mathbf{z}, \mathbf{x}_t) = \max_{\substack{\tilde{\mathbf{y}}, \tilde{\mathbf{x}}, \\ \|\tilde{\mathbf{x}} - \mathbf{x}_t\| \leq \Gamma, \\ \tilde{\mathbf{y}} \in Y_\Gamma(\mathbf{x}_t)}} \left\{ \min_{\substack{\Delta_t, \mathbf{s}_t^-, \mathbf{s}_t^+, \\ \mathbf{f}_t^{wc}, \beta_t^{wc}}} M^-\left(\mathbf{e}^\top \mathbf{s}_t^-\right) + M^+\left(\mathbf{e}^\top \mathbf{s}_t^+\right) \right. \tag{28k}$$

$$s.t. \quad \mathbf{G}\Delta_t + \mathbf{A}\mathbf{f}_t^{wc} - \mathbf{s}_t^+ + \mathbf{s}_t^-$$
$$= \mathbf{D}_t - \mathbf{E}\tilde{\mathbf{y}} - \mathbf{G}\mathbf{p}_t, \tag{28l}$$

$$\mathbf{f}_t^{wc} = \mathbf{B}\beta_t^{wc}, \tag{28m}$$

$$\mathbf{r}_t^{dn} \leq \Delta_t \leq \mathbf{r}_t^{up}, \tag{28n}$$

$$-\overline{\mathbf{F}} \leq \mathbf{f}_t^{wc} \leq \overline{\mathbf{F}}, \tag{28o}$$

$$0 \leq \mathbf{s}_t^- \leq \mathbf{D}_t, \tag{28p}$$

$$\left. 0 \leq \mathbf{s}_t^+ \leq \mathbf{E}\tilde{\mathbf{y}} \right\} \tag{28q}$$

*The online energy application problem (28) is formulated as a trilevel optimization problem (min-max-min). The objective function (28a) comprises production costs, up- and down-reserve costs. The constraint (28b) represents the nodal power balance in a dc power flow model. Expression (28c) represents Kirchhoff's second law. Limits for energy levels and up- and down-spinning reserves are imposed in the constraint (28d). The limits for the transmission line power flow, up- and down-spinning are imposed in expressions (28e) - (28g). Inter-period ramping limits are modeled by (28h) - (28i). Equation (15) ensures re-dispatch capability within a feasibility penalty cost $\alpha$ under the set $Y_\Gamma(\mathbf{x})$ for all plausible realizations within the DDCU set.*

*The second-level problem encompasses the outer maximization challenge as described in (28k). This problem aims to determine the weights $\theta = (\theta_1, \ldots, \theta_S)^\top$ for a convex combination $\tilde{\mathbf{x}}$, which is conditioned upon the new observation of contextual information $\mathbf{x}$. This means identifying a convex combination $\tilde{\mathbf{x}}$ of the covariate scenarios such that the distance $|\tilde{\mathbf{x}} - \mathbf{x}|$ is defined within the constraint of the robustness budget for contextual information $\Gamma$, i.e., $|\tilde{\mathbf{x}} - \mathbf{x}| \leq \Gamma$. This process ultimately defines the worst-case uncertainty realization $\tilde{\mathbf{h}}$ associated to $\tilde{\mathbf{x}}$ for a given $[\mathbf{p}_t, \mathbf{r}_t^{up}, \mathbf{r}_t^{dn}]$. The measure of worst case is given by the minimum energy imbalance cost function, which receives as inputs the values of the first- and second-level variables $[\mathbf{p}_t, \mathbf{r}_t^{up}, \mathbf{r}_t^{dn}]$ and $\tilde{\mathbf{h}}$. The third-level problem, i.e., the inner minimization problem (28k) - (28q), plays the role of the minimum power imbalance cost function. The objective of this problem is to find a re-dispatch solution that minimizes the total cost of the load shedding $\mathbf{s}_t^-$ and the REG spillage $\mathbf{s}_t^-$ introduced in the power balance constraint (28l). Constraints (28m) and (28o) are to first-level expressions (28c) and (28e), respectively. Expression (28n) imposes that the re-dispatch $\Delta_t$ must be between the down- and up-reserve. Limits for the load shedding and REG spillage are imposed in expressions (28p) - (28q).*